
%
%
%
%
%
%
\magnification=\magstephalf      
%
%
\vsize=7.5truein                 
\hsize=5.2truein                 
\newskip\stdskip                 
\stdskip=6pt plus3pt minus3pt    
\medskipamount=\stdskip          
\parindent=0pt                   
\parskip=\stdskip                
\abovedisplayskip=\stdskip       
\belowdisplayskip=\stdskip       
\mathsurround=0.75pt             
\overfullrule=0pt                
%
%
\def\ppar{\par\goodbreak\vskip 8pt plus 4pt minus 4pt}     
%
%
\def\stdspace{\hskip 0.75em plus 0.15em\ignorespaces}
\let\qua\stdspace 
%
%
%
%
%
%
%
\def\hexnumber#1{\ifcase#1 0\or 1\or 2\or 3\or 4\or 5\or 6\or 7\or 8\or
 9\or A\or B\or C\or D\or E\or F\fi}
%
%
\font\thirtnmsa=msam10 scaled 1315    
\font\tenmsa=msam10          \font\ninemsa=msam9
\font\sevenmsa=msam7         \font\sixmsa=msam6
\font\fivemsa=msam5
%
%
\newfam\msafam                  \textfont\msafam=\tenmsa
\scriptfont\msafam=\sevenmsa    \scriptscriptfont\msafam=\fivemsa
\edef\hexa{\hexnumber\msafam}        
\def\msa{\fam\msafam\tenmsa}         
%
%
\font\thirtnmsb=msbm10 scaled 1315   
\font\tenmsb=msbm10      \font\ninemsb=msbm9
\font\sevenmsb=msbm7     \font\sixmsb=msbm6
\font\fivemsb=msbm5
%
\newfam\msbfam                   \textfont\msbfam=\tenmsb       
\scriptfont\msbfam=\sevenmsb     \scriptscriptfont\msbfam=\fivemsb
\edef\hexb{\hexnumber\msbfam}    
\def\msb{\fam\msbfam\tenmsb}     
%
%
\font\thirtneufm=eufm10 scaled 1315   
\font\teneufm=eufm10                 \font\nineeufm=eufm9
\font\seveneufm=eufm7                \font\sixeufm=eufm6
\font\fiveeufm=eufm5
%
\newfam\eufmfam                    \textfont\eufmfam=\teneufm
\scriptfont\eufmfam=\seveneufm     \scriptscriptfont\eufmfam=\fiveeufm
\edef\hexf{\hexnumber\eufmfam}      
\def\frak{\fam\eufmfam\teneufm}     
%
%
%
\font\thirtnrm=cmr10 scaled 1315    
\font\ninerm=cmr9                   \font\sixrm=cmr6   
%
\font\thirtni=cmmi10 scaled 1315    
\font\ninei=cmmi9                   \font\sixi=cmmi6  
%
\font\thirtnsy=cmsy10 scaled 1315   
\font\ninesy=cmsy9                  \font\sixsy=cmsy6  
%
\font\thirtnbf=cmbx10 scaled 1315   
\font\ninebf=cmbx9                  \font\sixbf=cmbx6  
%
%
\font\thirtnex=cmex10 scaled 1315   
\font\nineex=cmex9                  
%
%
\font\thirtnit=cmti10 scaled 1315  
\font\nineit=cmti9                  
%
\font\thirtnsl=cmsl10 scaled 1315  
\font\ninesl=cmsl9                  
%
\font\thirtntt=cmtt10 scaled 1315  
\font\ninett=cmtt9                  
%
%
%
%
\def\small{%
%
%
\textfont0=\ninerm \scriptfont0=\sixrm \scriptscriptfont0=\fiverm
\def\rm{\fam0\ninerm}
%
%
\textfont1=\ninei \scriptfont1=\sixi \scriptscriptfont1=\fivei
%
%
\textfont2=\ninesy \scriptfont2=\sixsy \scriptscriptfont2=\fivesy
%
%
\textfont3=\nineex \scriptfont3=\nineex \scriptscriptfont3=\nineex
%
%
\textfont\bffam=\ninebf \scriptfont\bffam=\sixbf
\scriptscriptfont\bffam=\fivebf \def\bf{\fam\bffam\ninebf}%
%
%
\textfont\itfam=\nineit \def\it{\fam\itfam\nineit}%
\textfont\slfam=\ninesl \def\sl{\fam\slfam\ninesl}%
\textfont\ttfam=\ninett \def\tt{\fam\ttfam\ninett}%
%
%
%
\textfont\msafam=\ninemsa \scriptfont\msafam=\sixmsa
\scriptscriptfont\msafam=\fivemsa \def\msa{\fam\msafam\ninemsa}%
%
%
\textfont\msbfam=\ninemsb \scriptfont\msbfam=\sixmsb
\scriptscriptfont\msbfam=\fivemsb \def\msb{\fam\msbfam\ninemsb}%
%
%
\textfont\eufmfam=\nineeufm  \scriptfont\eufmfam=\sixeufm
\scriptscriptfont\eufmfam=\fiveeufm \def\frak{\fam\eufmfam\nineeufm}%
%
%
%
\normalbaselineskip=11pt%
\setbox\strutbox=\hbox{\vrule height8pt depth3pt width0pt}%
%
%
\normalbaselines\rm
%
%
\stdskip=4pt plus2pt minus2pt    
\medskipamount=\stdskip          
\parskip=\stdskip                
\abovedisplayskip=\stdskip       
\belowdisplayskip=\stdskip       
\def\ppar{\par\goodbreak\vskip 6pt plus 3pt minus 3pt}%
%
%
\def\section##1{\global\advance\sectionnumber by 1
\vskip-\lastskip\penalty-800\vskip 20pt plus10pt minus5pt 
\egroup{\bf\number\sectionnumber\quad##1}\bgroup\small         
\vskip 6pt plus3pt minus3pt
\nobreak\resultnumber=1}
}    
%
\def\beginsmall{\bgroup\small}
\let\endsmall\egroup
%
%
%
%
\def\large{%
\textfont0=\thirtnrm \scriptfont0=\ninerm \scriptscriptfont0=\sevenrm
\def\rm{\fam0\thirtnrm}%
\textfont1=\thirtni \scriptfont1=\ninei \scriptscriptfont1=\seveni
\textfont2=\thirtnsy \scriptfont2=\ninesy \scriptscriptfont2=\sevensy
\textfont3=\thirtnex \scriptfont3=\thirtnex \scriptscriptfont3=\thirtnex
\textfont\bffam=\thirtnbf \scriptfont\bffam=\ninebf
\scriptscriptfont\bffam=\sevenbf \def\bf{\fam\bffam\thirtnbf}%
\textfont\itfam=\thirtnit \def\it{\fam\itfam\thirtnit}%
\textfont\slfam=\thirtnsl \def\sl{\fam\slfam\thirtnsl}%
\textfont\ttfam=\thirtntt \def\tt{\fam\ttfam\thirtntt}%
\textfont\msafam=\thirtnmsa \scriptfont\msafam=\ninemsa
\scriptscriptfont\msafam=\sevenmsa \def\msa{\fam\msafam\thirtnmsa}%
\textfont\msbfam=\thirtnmsb \scriptfont\msbfam=\ninemsb
\scriptscriptfont\msbfam=\sevenmsb \def\msb{\fam\msbfam\thirtnmsb}%
\textfont\eufmfam=\thirtneufm  \scriptfont\eufmfam=\nineeufm
\scriptscriptfont\eufmfam=\seveneufm \def\frak{\fam\eufmfam\teneufm}%
\normalbaselineskip=16pt%
\setbox\strutbox=\hbox{\vrule height11.5pt depth4.5pt width0pt}%
\normalbaselines\rm}%
\let\Large\large   
%
\def\Bbb#1{{\msb#1}}

%

%
\mathchardef\plussquare="0\hexa01
\mathchardef\nge="3\hexb0B
\mathchardef\maltesecross="0\hexa7A
\mathchardef\del="0\hexf01
%
%
%
%
\font\sc=cmcsc10
%
%
%
%
\def\sqr#1#2{{\vcenter{\vbox{\hrule  height.#2truept
	\hbox{\vrule width.#2truept height#1truept 
	\kern#1truept \vrule width.#2truept}
	\hrule height.#2truept}}}}
\def\sq{\sqr55}    
%
%
%
%
\newcount\sectionnumber            
\newcount\resultnumber             
\sectionnumber=0\resultnumber=1    
%
%
%
\def\section#1{\global\advance\sectionnumber by 1
\xdef\nextkey{\number\sectionnumber}
\vskip-\lastskip\penalty-800\vskip 20pt plus10pt minus5pt 
{\large\bf\number\sectionnumber\quad#1}         
\vskip 8pt plus4pt minus4pt
\nobreak\resultnumber=1}      
%
%
%
%
%
\def\sh#1{\vskip-\lastskip\ppar{\bf #1}\par\nobreak\medskip}         
%
%
%
%

%
\def\proc#1{\xdef\nextkey{\number\sectionnumber.\number\resultnumber}%
\vskip-\lastskip\ppar\bf%
\noindent#1\ \number\sectionnumber.\number\resultnumber
\stdspace\sl\global\advance\resultnumber by 1\ignorespaces}
 
%
%
\def\qed{\hfill$\sq$\par\goodbreak\rm}   
%
%
%
%
%
%
%
%
\def\proclaim#1{\vskip-\lastskip\ppar\bf%
\noindent#1\stdspace\sl\ignorespaces} 
\let\endproclaim\endproc
%
%
%
%
\def\rk#1{\vskip-\lastskip\ppar{\bf #1}\stdspace\ignorespaces}                

%
%
%
%
%
%
\def\label{\xdef\nextkey{\number\sectionnumber.\number\resultnumber}%
\number\sectionnumber.\number\resultnumber
\global\advance\resultnumber by 1}
%
%
%
%
%
%
%
%
%
%
%
%
%
%
%
%
\newcount\refnumber              
\refnumber=1                     
\long\def\reflist#1\endreflist{%
\long\def\thereflist{#1}{\def\refkey##1##2\par{\xdef##1{\number\refnumber}%
\global\advance\refnumber by 1}%
\def\key##1##2\par{\expandafter\xdef%
\csname##1\endcsname{\number\refnumber}%
\global\advance\refnumber by 1}#1\par}}
\long\def\references{%
\penalty-800\vskip-\lastskip\vskip 15pt plus10pt minus5pt 
{\large\bf References}\ppar 
{\leftskip=25pt\frenchspacing    
\small\parskip=3pt plus2pt       
\def\refkey##1##2\par{\noindent  
\llap{[##1]\stdspace}\ignorespaces##2\par}         
\def\key##1##2\par{\noindent  
\llap{[\ref{##1}]\stdspace}\ignorespaces##2\par}  
\def\,{\thinspace}\thereflist\par}}
%
%
%
\newcount\footnotenumber         
\footnotenumber=1                
\def\fnote#1{\xdef\nextkey{\number\footnotenumber}%
{\small\ifnum\footnotenumber>9\parindent=14pt%
\else\parindent=10pt\fi\footnote{$^{\number\footnotenumber}$}%
{\hglue-5pt#1}\global\advance\footnotenumber by 1}}
%
%
%
%
%
%
%
\newcount\figurenumber          
\figurenumber=1                 
\def\caption#1{\xdef\nextkey{\number\figurenumber}%
\cl{\small Figure \number\figurenumber: #1}%
\global\advance\figurenumber by 1}
\def\figurelabel{\xdef\nextkey{\number\figurenumber}%
\cl{\small Figure \number\figurenumber}%
\global\advance\figurenumber by 1}
\long\def\figure#1\endfigure{{\xdef\nextkey{\number\figurenumber}%
\let\captiontext\relax\def\caption##1{\xdef\captiontext{##1}}%
\midinsert\cl{\ignorespaces#1\unskip\unskip\unskip\unskip}\vglue6pt\cl{\small 
Figure \number\figurenumber\ifx\captiontext\relax\else: \captiontext
\fi}\endinsert\global\advance\figurenumber by 1}}
%
%
%
%
%
%
%
\def\nextkey{??}   
%
\def\key#1{\expandafter\xdef\csname #1\endcsname{\nextkey}}
\def\ref#1{\expandafter\ifx\csname #1\endcsname\relax
\immediate\write16{Reference {#1} undefined}??\else
\csname #1\endcsname\fi}
%
%
%
%
%
%
%
\newread\gtinfile
\newwrite\gtreffile
\def\useforwardrefs{
\openin\gtinfile\jobname.ref
\ifeof\gtinfile
\closein\gtinfile
\immediate\write16{No file \jobname.ref}
\else
\closein\gtinfile
\input \jobname.ref
\fi
\immediate\openout\gtreffile \jobname.ref
%
%
\def\key##1{{\def\\{\noexpand}%
\expandafter\xdef\csname ##1\endcsname{\nextkey}%
\immediate\write\gtreffile{\\\expandafter\\\def\\\csname ##1\\\endcsname%
{\nextkey}}}}
%
%
\long\def\reflist##1\endreflist{%
\long\def\thereflist{##1}{\def\refkey####1####2\par{\xdef####1{%
\number\refnumber}{\def\\{\noexpand}\immediate\write\gtreffile
{\\\def\\####1{\number\refnumber}}}\global\advance\refnumber by 1}%
\def\key####1####2\par{\expandafter\xdef%
\csname####1\endcsname{\number\refnumber}%
{\def\\{\noexpand}\immediate\write\gtreffile
{\\\expandafter\\\def\\\csname ####1\\\endcsname{\number\refnumber}}}
\global\advance\refnumber by 1}##1\par}}
\long\def\biblio##1\endbiblio{\reflist##1\endreflist\references}%
%
%
\def\numkey##1{{\def\\{\noexpand}%
\xdef##1{\number\sectionnumber.\number\resultnumber}
\immediate\write\gtreffile{\\\def\\##1%
{\number\sectionnumber.\number\resultnumber}}}}
\def\seckey##1{{\def\\{\noexpand}\xdef##1{\number\sectionnumber}
\immediate\write\gtreffile{\\\def\\##1{\number\sectionnumber}}}}
\def\figkey##1{\xdef##1{\number\figurenumber}%
{\def\\{\noexpand}\immediate\write\gtreffile%
{\\\def\\##1{\number\figurenumber}}}
\number\figurenumber\global\advance\figurenumber by 1}
}   
%
%
%
%
\def\figkey#1{\xdef#1{\number\figurenumber}%
\number\figurenumber\global\advance\figurenumber by 1}
\def\fig#1#2\endfig{%
\midinsert\cl{#2}\vglue6pt\cl{\small Figure #1}\endinsert}
\def\newfig{\number\figurenumber\global\advance\figurenumber by 1}
\def\numkey#1{\xdef#1{\number\sectionnumber.\number\resultnumber}}
\def\seckey#1{\xdef#1{\number\sectionnumber}}
%
%
%
%
%
%
%
%
%
\def\verb{\catcode`\"=\active}       
\def\brev{\catcode`\"=12}            
\brev                                
\verb                                
{\obeyspaces\gdef {\ }}              
{\catcode`\`=\active\gdef`{\relax\lq}}
\def"{%
\begingroup\baselineskip=12pt\def\par{\leavevmode\endgraf}%
\tt\obeylines\obeyspaces\parskip=0pt\parindent=0pt%
\catcode`\$=12\catcode`\&=12\catcode`\^=12\catcode`\#=12%
\catcode`\_=12\catcode`\~=12%
\catcode`\{=12\catcode`\}=12\catcode`\%=12\catcode`\\=12%
\catcode`\`=\active\let"\endgroup}
\brev      
%
%
%
%
%
%
\def\items{\par\leftskip = 25pt}           
\def\enditems{\par\leftskip = 0pt}         
\def\item#1{\par\leavevmode\llap{#1\stdspace}%
\ignorespaces}                             
%
%

%
%
\def\co{\colon\thinspace}    
\def\np{\vfil\eject}         
\def\nl{\hfil\break}         
\def\cl{\centerline}         
\def\gt{{\mathsurround=0pt\it $\cal G\mskip-2mu$eometry \&\ 
$\cal T\!\!$opology}}        
\def\agt{{\mathsurround=0pt\it$\cal A\mskip-.7mu$lgebraic \&\ 
$\cal G\mskip-2mu$eometric $\cal T\!\!$opology}}  
%
%
%

%
%
%
%
%
\def\title#1{\def\thetitle{#1}}

\def\author#1{\edef\previousauthors{\theauthors}
 \ifx\theauthors\relax\def\theauthors{#1}\else
 \def\theauthors{\previousauthors\par#1}\fi}

%
\def\address#1{\edef\previousaddresses{\theaddress}
 \ifx\theaddress\relax\def\theaddress{#1}\else
 \def\theaddress{\previousaddresses\par\vskip 2pt\par#1}\fi}
\def\secondaddress#1{\edef\previousaddresses{\theaddress}
 \ifx\theaddress\relax\def\theaddress{#1}\else
 \def\theaddress{\previousaddresses\par{\rm and}\par#1}\fi}   

\def\email#1{\edef\previousemails{\theemail}
 \ifx\theemail\relax\def\theemail{#1}\else
 \def\theemail{\previousemails\hskip 0.75em\relax#1}\fi}
\def\secondemail#1{\edef\previousemails{\theemail}
 \ifx\theemail\relax\def\theemail{#1}\else
 \def\theemail{\previousemails\hskip 0.75em{\rm and}\hskip 0.75em
 \relax#1}\fi}
\def\url#1{\edef\previousurls{\theurl}
 \ifx\theurl\relax\def\theurl{#1}\else
 \def\theurl{\previousurls\hskip 0.75em\relax#1}\fi}
\def\secondurl#1{\edef\previousurls{\theurl}
 \ifx\theurl\relax\def\theurl{#1}\else
 \def\theurl{\previousurls\hskip 0.75em{\rm and}\hskip 0.75em
 \relax#1}\fi}
\long\def\abstract#1\endabstract{\long\def\theabstract{#1}}
\def\primaryclass#1{\def\theprimaryclass{#1}}
\let\subjclass\primaryclass                        
\def\secondaryclass#1{\def\thesecondaryclass{#1}}
\def\keywords#1{\def\thekeywords{#1}}
%
%
\let\\\par\let\thetitle\relax\let\theshorttitle\relax
\let\theauthors\relax\let\theshortauthors\relax
\let\theaddress\relax\let\theshortaddress\relax
\let\theemail\relax\let\theurl\relax
\let\theabstract\relax\let\theprimaryclass\relax
\let\thesecondaryclass\relax\let\thekeywords\relax
%
%
%
%
\long\def\maketitlepage{    

\vglue 0.2truein   

%
{\parskip=0pt\leftskip 0pt plus 1fil\def\\{\par\smallskip}{\large
\bf\thetitle}\par\medskip}   

\vglue 0.15truein 

%
{\parskip=0pt\leftskip 0pt plus 1fil\def\\{\par}{\sc\theauthors}
\par\medskip}%
 
\vglue 0.1truein 

%
{\small\parskip=0pt
{\leftskip 0pt plus 1fil\def\\{\par}{\sl\theaddress}\par}
\ifx\theemail\relax\else  
\vglue 5pt \def\\{\stdspace{\rm and}\stdspace} 
\cl{Email:\stdspace\tt\theemail}\fi
\ifx\theurl\relax\else    
\vglue 5pt \def\\{\stdspace{\rm and}\stdspace} 
\cl{URL:\stdspace\tt\theurl}\fi\par}

\vglue 7pt 

{\bf Abstract}

\vglue 5pt

\theabstract

\vglue 7pt 

{\bf AMS Classification numbers}\quad Primary:\quad \theprimaryclass\par

Secondary:\quad \thesecondaryclass

\vglue 5pt 

{\bf Keywords:}\quad \thekeywords

\np  

}    
%
%
\long\def\makeshorttitle{    


%
{\parskip=0pt\leftskip 0pt plus 1fil\def\\{\par\smallskip}{\large
\bf\thetitle}\par\medskip}   

\vglue 0.05truein 

%
{\parskip=0pt\leftskip 0pt plus 1fil\def\\{\par}{\sc\theauthors}
\par\medskip}%
 
\vglue 0.03truein 

%
{\small\parskip=0pt
{\leftskip 0pt plus 1fil\def\\{\par}{\sl\ifx\theshortaddress\relax
\theaddress\else\theshortaddress\fi}\par}
\ifx\theemail\relax\else  
\vglue 5pt \def\\{\stdspace{\rm and}\stdspace} 
\cl{Email:\stdspace\tt\theemail}\fi
\ifx\theurl\relax\else    
\vglue 5pt \def\\{\stdspace{\rm and}\stdspace} 
\cl{URL:\stdspace\tt\theurl}\fi\par}

\vglue 10pt 


{\small\leftskip 25pt\rightskip 25pt{\bf Abstract}\stdspace\theabstract

{\bf AMS Classification}\stdspace\theprimaryclass
\ifx\thesecondaryclass\relax\else; \thesecondaryclass\fi\par
{\bf Keywords}\stdspace \thekeywords\par}
\vglue 7pt
}    
\let\maketitle\makeshorttitle        
%
%

\def\volumenumber#1{\def\thevolumenumber{#1}}
\def\volumeyear#1{\def\thevolumeyear{#1}}
\def\pagenumbers#1#2{\def\startpage{#1}\def\finishpage{#2}}
\def\published#1{\def\publishdate{#1}}
\def\received#1{\def\receiveddate{#1}}
\def\revised#1{\def\reviseddate{#1}}
\let\reviseddate\relax
\volumenumber{X}
\volumeyear{20XX}
\pagenumbers{1}{XXX}
\published{XX Xxxember 20XX}

\long\def\makeagttitle{   
\agt\hfill      
\hbox to 60truept{\vbox to 0pt{\vglue -14truept{\bf [Logo here]}\vss}\hss}
\break
{\small Volume \thevolumenumber\ (\thevolumeyear)
\startpage--\finishpage\nl
Published: \publishdate}

\vglue .2truein

{\parskip=0pt\leftskip 0pt plus 1fil\def\\{\par\smallskip}{\large
\bf\thetitle}\par\medskip}   
\vglue 0.05truein 

%
{\parskip=0pt\leftskip 0pt plus 1fil\def\\{\par}{\sc\theauthors}
\par\medskip}%
 
\vglue 0.03truein 


{\small\leftskip 25truept\rightskip 25truept{\bf Abstract}\stdspace\theabstract

{\bf AMS Classification}\stdspace\theprimaryclass
\ifx\thesecondaryclass\relax\else; \thesecondaryclass\fi\par
{\bf Keywords}\stdspace \thekeywords\par}\vglue 7truept

}   


\def\Addresses{\bigskip
{\small \parskip 0pt \leftskip 0pt \rightskip 0pt plus 1fil \def\\{\par}
\sl\theaddress\par\medskip \rm Email:\stdspace\tt\theemail\par
\ifx\theurl\relax\else\smallskip \rm URL:\stdspace\tt\theurl\par\fi}}

\def\agtart{
\hoffset 14truemm
\voffset 31truemm
\font\phead=cmsl9 scaled 950
\font\pnum=cmbx10 scaled 913
\font\pfoot=cmsl9 scaled 950
\headline{\vbox to 0pt{\vskip -4.5mm\line{\small\phead\ifnum
\count0=\startpage ISSN numbers are printed here
\hfill {\pnum\folio}\else\ifodd\count0\def\\{ }%
\ifx\theshorttitle\relax\thetitle\else\theshorttitle\fi\hfill{\pnum\folio}
\else\def\\{ and }{\pnum\folio}\hfill\ifx\theshortauthors\relax\theauthors
\else\theshortauthors\fi\fi\fi}\vss}}
\footline{\vbox to 0pt{\vglue 0mm\line{\small\pfoot\ifnum\count0=\startpage
Copyright declaration is printed here\hfill\else
\agt, Volume \thevolumenumber\ (\thevolumeyear)\hfill\fi}\vss}}
\let\maketitle\makeagttitle\let\makeshorttitle\makeagttitle}

%
%
%
%
%
%
%
\def\newsymbol#1#2#3#4#5{%
\ifodd#2\let\hexno\hexa\else\let\hexno\hexb\fi
\mathchardef#1="#3\hexno#4#5}
%
%
\def\mathsurroundoff{\mathsurround=0pt}
\def\mathhexbox#1#2#3{\relax
\ifmmode\mathpalette{}{\mathsurroundoff\mathchar"#1#2#3}%
\else\leavevmode\hbox{$\mathsurroundoff\mathchar"#1#2#3$}\fi}
%
%
\mathchardef\dashbar"0\hexa39

%
%
\newsymbol\boxdot 1200
\newsymbol\boxplus 1201
\newsymbol\boxtimes 1202
\newsymbol\square 1003
\newsymbol\blacksquare 1004
\newsymbol\centerdot 1205
\newsymbol\lozenge 1006
\newsymbol\blacklozenge 1007
\newsymbol\circlearrowright 1308
\newsymbol\circlearrowleft 1309
\let\rightleftharpoons\undefined
\newsymbol\rightleftharpoons 130A
\newsymbol\leftrightharpoons 130B
\newsymbol\boxminus 120C
\newsymbol\Vdash 130D
\newsymbol\Vvdash 130E
\newsymbol\vDash 130F
\newsymbol\twoheadrightarrow 1310
\newsymbol\twoheadleftarrow 1311
\newsymbol\leftleftarrows 1312
\newsymbol\rightrightarrows 1313
\newsymbol\upuparrows 1314
\newsymbol\downdownarrows 1315
\newsymbol\upharpoonright 1316
 
\newsymbol\downharpoonright 1317
\newsymbol\upharpoonleft 1318
\newsymbol\downharpoonleft 1319
\newsymbol\rightarrowtail 131A
\newsymbol\leftarrowtail 131B
\newsymbol\leftrightarrows 131C
\newsymbol\rightleftarrows 131D
\newsymbol\Lsh 131E
\newsymbol\Rsh 131F
\newsymbol\rightsquigarrow 1320
\newsymbol\leftrightsquigarrow 1321
\newsymbol\looparrowleft 1322
\newsymbol\looparrowright 1323
\newsymbol\circeq 1324
\newsymbol\succsim 1325
\newsymbol\gtrsim 1326
\newsymbol\gtrapprox 1327
\newsymbol\multimap 1328
\newsymbol\therefore 1329
\newsymbol\because 132A
\newsymbol\doteqdot 132B
 
\newsymbol\triangleq 132C
\newsymbol\precsim 132D
\newsymbol\lesssim 132E
\newsymbol\lessapprox 132F
\newsymbol\eqslantless 1330
\newsymbol\eqslantgtr 1331
\newsymbol\curlyeqprec 1332
\newsymbol\curlyeqsucc 1333
\newsymbol\preccurlyeq 1334
\newsymbol\leqq 1335
\newsymbol\leqslant 1336
\newsymbol\lessgtr 1337
\newsymbol\backprime 1038
\newsymbol\risingdotseq 133A
\newsymbol\fallingdotseq 133B
\newsymbol\succcurlyeq 133C
\newsymbol\geqq 133D
\newsymbol\geqslant 133E
\newsymbol\gtrless 133F
\newsymbol\sqsubset 1340
\newsymbol\sqsupset 1341
\newsymbol\vartriangleright 1342
\newsymbol\vartriangleleft 1343
\newsymbol\trianglerighteq 1344
\newsymbol\trianglelefteq 1345
\newsymbol\bigstar 1046
\newsymbol\between 1347
\newsymbol\blacktriangledown 1048
\newsymbol\blacktriangleright 1349
\newsymbol\blacktriangleleft 134A
\newsymbol\vartriangle 134D
\newsymbol\blacktriangle 104E
\newsymbol\triangledown 104F
\newsymbol\eqcirc 1350
\newsymbol\lesseqgtr 1351
\newsymbol\gtreqless 1352
\newsymbol\lesseqqgtr 1353
\newsymbol\gtreqqless 1354
\newsymbol\Rrightarrow 1356
\newsymbol\Lleftarrow 1357
\newsymbol\veebar 1259
\newsymbol\barwedge 125A
\newsymbol\doublebarwedge 125B
\let\angle\undefined
\newsymbol\angle 105C
\newsymbol\measuredangle 105D
\newsymbol\sphericalangle 105E
\newsymbol\varpropto 135F
\newsymbol\smallsmile 1360
\newsymbol\smallfrown 1361
\newsymbol\Subset 1362
\newsymbol\Supset 1363
\newsymbol\Cup 1264
 
\newsymbol\Cap 1265
 
\newsymbol\curlywedge 1266
\newsymbol\curlyvee 1267
\newsymbol\leftthreetimes 1268
\newsymbol\rightthreetimes 1269
\newsymbol\subseteqq 136A
\newsymbol\supseteqq 136B
\newsymbol\bumpeq 136C
\newsymbol\Bumpeq 136D
\newsymbol\lll 136E
 
\newsymbol\ggg 136F
 
\newsymbol\circledS 1073
\newsymbol\pitchfork 1374
\newsymbol\dotplus 1275
\newsymbol\backsim 1376
\newsymbol\backsimeq 1377
\newsymbol\complement 107B
\newsymbol\intercal 127C
\newsymbol\circledcirc 127D
\newsymbol\circledast 127E
\newsymbol\circleddash 127F
\newsymbol\lvertneqq 2300
\newsymbol\gvertneqq 2301
\newsymbol\nleq 2302
\newsymbol\ngeq 2303
\newsymbol\nless 2304
\newsymbol\ngtr 2305
\newsymbol\nprec 2306
\newsymbol\nsucc 2307
\newsymbol\lneqq 2308
\newsymbol\gneqq 2309
\newsymbol\nleqslant 230A
\newsymbol\ngeqslant 230B
\newsymbol\lneq 230C
\newsymbol\gneq 230D
\newsymbol\npreceq 230E
\newsymbol\nsucceq 230F
\newsymbol\precnsim 2310
\newsymbol\succnsim 2311
\newsymbol\lnsim 2312
\newsymbol\gnsim 2313
\newsymbol\nleqq 2314
\newsymbol\ngeqq 2315
\newsymbol\precneqq 2316
\newsymbol\succneqq 2317
\newsymbol\precnapprox 2318
\newsymbol\succnapprox 2319
\newsymbol\lnapprox 231A
\newsymbol\gnapprox 231B
\newsymbol\nsim 231C
\newsymbol\ncong 231D
\newsymbol\diagup 201E
\newsymbol\diagdown 201F
\newsymbol\varsubsetneq 2320
\newsymbol\varsupsetneq 2321
\newsymbol\nsubseteqq 2322
\newsymbol\nsupseteqq 2323
\newsymbol\subsetneqq 2324
\newsymbol\supsetneqq 2325
\newsymbol\varsubsetneqq 2326
\newsymbol\varsupsetneqq 2327
\newsymbol\subsetneq 2328
\newsymbol\supsetneq 2329
\newsymbol\nsubseteq 232A
\newsymbol\nsupseteq 232B
\newsymbol\nparallel 232C
\newsymbol\nmid 232D
\newsymbol\nshortmid 232E
\newsymbol\nshortparallel 232F
\newsymbol\nvdash 2330
\newsymbol\nVdash 2331
\newsymbol\nvDash 2332
\newsymbol\nVDash 2333
\newsymbol\ntrianglerighteq 2334
\newsymbol\ntrianglelefteq 2335
\newsymbol\ntriangleleft 2336
\newsymbol\ntriangleright 2337
\newsymbol\nleftarrow 2338
\newsymbol\nrightarrow 2339
\newsymbol\nLeftarrow 233A
\newsymbol\nRightarrow 233B
\newsymbol\nLeftrightarrow 233C
\newsymbol\nleftrightarrow 233D
\newsymbol\divideontimes 223E
\newsymbol\varnothing 203F
\newsymbol\nexists 2040
\newsymbol\Finv 2060
\newsymbol\Game 2061
\newsymbol\mho 2066
\newsymbol\eth 2067
\newsymbol\eqsim 2368
\newsymbol\beth 2069
\newsymbol\gimel 206A
\newsymbol\daleth 206B
\newsymbol\lessdot 236C
\newsymbol\gtrdot 236D
\newsymbol\ltimes 226E
\newsymbol\rtimes 226F
\newsymbol\shortmid 2370
\newsymbol\shortparallel 2371
\newsymbol\smallsetminus 2272
\newsymbol\thicksim 2373
\newsymbol\thickapprox 2374
\newsymbol\approxeq 2375
\newsymbol\succapprox 2376
\newsymbol\precapprox 2377
\newsymbol\curvearrowleft 2378
\newsymbol\curvearrowright 2379
\newsymbol\digamma 207A
\newsymbol\varkappa 207B
\newsymbol\Bbbk 207C
\newsymbol\hslash 207D
\let\hbar\undefined
\newsymbol\hbar 207E
\newsymbol\backepsilon 237F

\input amstex
\input epsf
%
%
\def\relabelbox{%
  \hbox\bgroup%
}%
\def\endrelabelbox{%
\egroup%
}%
\def\relabel #1#2 {%
  \special{ps:/a {} def}%
  \smash{\rlap{#2}}%
}%
\def\adjustrelabel <#1,#2> #3#4 {%
  \special{ps:/a {} def}%
  \smash{\rlap{\kern #1 \raise #2\hbox{#4}}}%
}%
\def\extralabel <#1,#2> #3 {\smash{\rlap{\kern #1 \raise #2\hbox{#3}}}}%

\chardef\newinsCatAt\the\catcode `\@
\catcode `\@=11
%
%
%
\newskip\insertskipamount\newskip\inserthardskipamount
\insertskipamount 12pt plus2pt     
\inserthardskipamount 4pt          
\def\insertskip{\vskip\insertskipamount}
%
%
\newskip\LastSkip
\def\SaveLastSkip{\LastSkip\lastskip}
\def\RestoreLastSkip{\nobreak\vskip-\LastSkip\vskip\LastSkip}
%
%
\newcount\SplitTest
\def\SetSplitTest{\SplitTest\insertpenalties
  \insert\topins{\floatingpenalty1}%
  \advance\SplitTest-\insertpenalties}
%
%
\def\midinsert{\par
 \SaveLastSkip\penalty-150\SetSplitTest\RestoreLastSkip
 \ifnum\SplitTest=-1
  \@midfalse\p@gefalse\else\@midtrue\fi\@ins}
\def\@ins{\par\begingroup\setbox\z@\vbox\bgroup%
  \vglue\inserthardskipamount}
\def\endinsert{\egroup 
  \if@mid \dimen@\ht\z@ \advance\dimen@\dp\z@
    \advance\dimen@\insertskipamount
    \advance\dimen@\pagetotal\advance\dimen@-\pageshrink
    \ifdim\dimen@>\pagegoal\@midfalse\p@gefalse\fi\fi
  \if@mid%
    \ifdim\lastskip<\insertskipamount\removelastskip\insertskip\fi
    \nointerlineskip\box\z@\penalty-200\insertskip
  \else%
    \SaveLastSkip
    \insert\topins{\penalty100 
    \splittopskip\z@skip
    \splitmaxdepth\maxdimen \floatingpenalty\z@
    \ifp@ge \dimen@\dp\z@
    \vbox to\vsize{\unvbox\z@\kern-\dimen@}
    \else \box\z@\nobreak\insertskip\fi}
    \RestoreLastSkip
   \fi\endgroup}
%
\catcode `\@=\newinsCatAt

\catcode`\@=12    


\def\ifplaintex{\expandafter\ifx\csname documentclass\endcsname\relax}


\ifplaintex 
\hoffset 14truemm
\voffset 31truemm
\else
\headsep 23pt
\footskip 35pt
\hoffset -4truemm
\voffset 12.5truemm
\fi

\expandafter\ifx\csname beginpicture\endcsname\relax
\expandafter\ifx\csname documentclass\endcsname\relax
\input pictex \else\font\fiverm=cmr5
\input prepictex \input pictex \input postpictex \fi\fi

\def\gt{{\mathsurround=0pt\it $\cal G\mskip-2mu$eometry \&\ 
$\cal T\!\!$opology}}        

\def\gtp{{\mathsurround=0pt\it $\cal G\mskip-2mu$eometry \&\ 
$\cal T\!\!$opology $\cal P\!$ublications}}  


\def\lognumber#1{\def\thelognumber{#1}}
\def\volumenumber#1{\def\thevolumenumber{#1}}
\def\papernumber#1{\def\thepapernumber{#1}}
\def\volumeyear#1{\def\thevolumeyear{#1}}

\def\pagenumbers#1#2{\def\startpage{#1}\def\finishpage{#2}}
\def\published#1{\def\publishdate{#1}}
\def\proposed#1{\def\theproposer{#1}}
\def\seconded#1{\def\theseconders{#1}}
\def\received#1{\def\receiveddate{#1}}
\def\revised#1{\def\reviseddate{#1}}
\def\accepted#1{\def\accepteddate{#1}}

\def\asciiaddress#1{\def\theasciiaddress{#1}}
\def\asciiemail#1{\def\theasciiemail{#1}}


\let\\\par\let\thelognumber\relax
\let\thevolumenumber\relax\let\thepapernumber\relax
\let\thevolumeyear\relax\let\thesamplenumber\relax\let\startpage\relax
\let\finishpage\relax\let\publishdate\relax\let\receiveddate\relax
\let\reviseddate\relax\let\accepteddate\relax\let\theasciititle\relax
\let\theasciiauthors\relax\let\theasciiaddress\relax
\let\theasciiabstract\relax
\let\theasciiemail\relax\let\theshortauthors\relax\let\theshorttitle\relax

\long\def\maketitlep{   

\count0=\startpage

\gt\hfill      
\beginpicture
\setcoordinatesystem units <0.33truein, 0.33truein> point at 2.2 0.9
\setplotsymbol ({$\cal G$})
\plotsymbolspacing=9truept
\circulararc 315 degrees from 0 1 center at 0 0
\setplotsymbol ({$\cal T$})
\circulararc 315 degrees from 1 -1 center at 1 0
\endpicture
%
\break
{\small\ifx\thesamplenumber\relax 
Volume \else Sample
\fi\thevolumenumber\ (\thevolumeyear)
\startpage--\finishpage\nl
Published: \publishdate}
\vglue 0.5truein plus 0.4fil minus 0.1truein

{\parskip=0pt\leftskip 0pt plus 1fil\def\\{\par\smallskip}{\ifplaintex\large
\else\Large\fi\bf\thetitle}\par\medskip}   

\vglue 0pt plus 0.1fil 

{\parskip=0pt\leftskip 0pt plus 1fil\def\\{\par}{\sc\theauthors}
\par\medskip}

\vglue 0pt plus 0.1fil 

{\small\parskip=0pt\let\newline\\
{\leftskip 0pt plus 1fil\def\\{\par}{\sl\theaddress}\par}
\expandafter\ifx\theemail\relax    
\relax\else\vglue 5pt plus 0.02fil minus 2pt\def\\{\stdspace{\rm 
and}\stdspace} 
\cl{Email:\stdspace\tt\theemail}\fi
\ifx\theurl\relax                  
\relax\else\vglue 5pt plus 0.02fil minus 2pt\def\\{\stdspace{\rm 
and}\stdspace}
\cl{URL:\stdspace\tt\theurl}\fi\par}

\vglue 7pt plus 0.3fil minus 3pt

{\bf Abstract}
\vglue 5pt plus 0.1fil minus 2pt

\theabstract

\vglue 7pt plus 0.3fil minus 3pt

{\bf AMS Classification numbers}\quad Primary:\quad \theprimaryclass

Secondary:\quad \thesecondaryclass

\vglue 5pt plus 0.3fil minus 2pt

{\bf Keywords}\quad \thekeywords

\vglue 10pt plus 0.5fil minus 5pt

{\small  Proposed: \theproposer\hfill Received: \receiveddate\nl
Seconded: \theseconders\hfill 
\ifx\reviseddate\relax                         
Accepted: \accepteddate                        
\else
Revised: \reviseddate                          
\fi}
\eject
}       

\let\maketitlepage\maketitlep
\let\maketitle\maketitlepage


\font\phead=cmsl9 scaled 950
\font\lhead=cmsl9 scaled 1050
\font\pnum=cmbx10 scaled 913
\font\lnum=cmbx10 
\font\pfoot=cmsl9 scaled 950
\font\lfoot=cmsl9 scaled 1050
\ifplaintex
\headline{\vbox to 0pt{\vskip -4.5mm\line{\small\phead\ifnum
\count0=\startpage ISSN 1364-0380 (on line)
1465-3060 (printed) \hfill {\pnum\folio}\else\ifodd\count0\def\\{ }%
\ifx\theshorttitle\relax\thetitle\else\theshorttitle\fi\hfill{\pnum\folio}
\else\def\\{ and }{\pnum\folio}\hfill\ifx\theshortauthors\relax\theauthors
\else\theshortauthors\fi\fi\fi}\vss}}
\footline{\vbox to 0pt{\vglue 0mm\line{\small\pfoot\ifnum\count0=\startpage
\copyright\ \gtp\hfill\else
\gt, Volume \thevolumenumber\ (\thevolumeyear)\hfill\fi}\vss
}}
\else
\makeatletter
\def\@oddhead{{\small\lhead\ifnum\count0=\startpage ISSN 1364-0380 (on line)
1465-3060 (printed) \hfill {\lnum\number\count0}\else\ifodd\count0
\def\\{ }\ifx\theshorttitle\relax \thetitle \else\theshorttitle\fi\hfill
{\lnum\number\count0}\else\def\\{ and }{\lnum\number\count0}
\hfill\ifx\theshortauthors\relax 
\theauthors\else\theshortauthors\fi\fi\fi}}\def\@evenhead{\@oddhead}
\def\@oddfoot{\small\lfoot\ifnum\count0=\startpage\copyright\ \gtp\hfill\else
\gt, Volume \thevolumenumber\ (\thevolumeyear)\hfill\fi}
\def\@evenfoot{\@oddfoot}
\makeatother
\fi


\newwrite\gtoutfile
\long\gdef\makeheadfile{  
{\def\\{, }\def\s{ }
\immediate\openout\gtoutfile head.xxx
\immediate\write\gtoutfile{To: math@arxiv.org}
\immediate\write\gtoutfile{Subject: put or rep NNNNN:pppp}
\immediate\write\gtoutfile{--text follows this line--}
\immediate\write\gtoutfile{Proxy-for: \ifx\theasciiauthors\relax
\theauthors\else\theasciiauthors\fi\s<\ifx\theasciiemail\relax\theemail\else\theasciiemail\fi>}
\immediate\write\gtoutfile{\noexpand\\}
\immediate\write\gtoutfile{Authors: \ifx\theasciiauthors\relax
\theauthors\else\theasciiauthors\fi}
\immediate\write\gtoutfile{Title: \ifx\theasciititle\relax
\thetitle\else\theasciititle\fi}
\immediate\write\gtoutfile{Subj-class: GT or SG or MG etc}
\immediate\write\gtoutfile{MSC-class: \theprimaryclass\ifx\thesecondaryclass\relax\else, \thesecondaryclass\fi}
\immediate\write\gtoutfile{Journal-ref: Geom. Topol. \thevolumenumber
(\thevolumeyear) \startpage-\finishpage}
\immediate\write\gtoutfile{Comments: Published by Geometry and Topology at}
\immediate\write\gtoutfile{\s\s http://www.maths.warwick.ac.uk/gt/GTVol\thevolumenumber/paper\thepapernumber.abs.html}
\immediate\write\gtoutfile{\noexpand\\}
\immediate\write\gtoutfile{}
\ifx\theasciiabstract\relax
\immediate\write\gtoutfile{\theabstract}\else
\immediate\write\gtoutfile{\theasciiabstract}\fi
\immediate\write\gtoutfile{}
\immediate\write\gtoutfile{\noexpand\\}
\immediate\write\gtoutfile{}
\immediate\closeout\gtoutfile}}  

\def\maketitlepage{\maketitlep\makeheadfile}
\let\maketitle\maketitlepage


\def\ifplaintex{\expandafter\ifx\csname documentclass\endcsname\relax}


\ifplaintex 
\hoffset 14truemm
\voffset 31truemm
\else
\headsep 23pt
\footskip 35pt
\hoffset -4truemm
\voffset 12.5truemm
\fi

\expandafter\ifx\csname beginpicture\endcsname\relax
\expandafter\ifx\csname documentclass\endcsname\relax
\input pictex \else\font\fiverm=cmr5
\input prepictex \input pictex \input postpictex \fi\fi

\def\gt{{\mathsurround=0pt\it $\cal G\mskip-2mu$eometry \&\ 
$\cal T\!\!$opology}}        

\def\gtp{{\mathsurround=0pt\it $\cal G\mskip-2mu$eometry \&\ 
$\cal T\!\!$opology $\cal P\!$ublications}}  


\def\lognumber#1{\def\thelognumber{#1}}
\def\volumenumber#1{\def\thevolumenumber{#1}}
\def\papernumber#1{\def\thepapernumber{#1}}
\def\volumeyear#1{\def\thevolumeyear{#1}}

\def\pagenumbers#1#2{\def\startpage{#1}\def\finishpage{#2}}
\def\published#1{\def\publishdate{#1}}
\def\proposed#1{\def\theproposer{#1}}
\def\seconded#1{\def\theseconders{#1}}
\def\received#1{\def\receiveddate{#1}}
\def\revised#1{\def\reviseddate{#1}}
\def\accepted#1{\def\accepteddate{#1}}

\def\asciiaddress#1{\def\theasciiaddress{#1}}
\def\asciiemail#1{\def\theasciiemail{#1}}


\let\\\par\let\thelognumber\relax
\let\thevolumenumber\relax\let\thepapernumber\relax
\let\thevolumeyear\relax\let\thesamplenumber\relax\let\startpage\relax
\let\finishpage\relax\let\publishdate\relax\let\receiveddate\relax
\let\reviseddate\relax\let\accepteddate\relax\let\theasciititle\relax
\let\theasciiauthors\relax\let\theasciiaddress\relax
\let\theasciiabstract\relax
\let\theasciiemail\relax\let\theshortauthors\relax\let\theshorttitle\relax

\long\def\maketitlep{   

\count0=\startpage

\gt\hfill      
\beginpicture
\setcoordinatesystem units <0.33truein, 0.33truein> point at 2.2 0.9
\setplotsymbol ({$\cal G$})
\plotsymbolspacing=9truept
\circulararc 315 degrees from 0 1 center at 0 0
\setplotsymbol ({$\cal T$})
\circulararc 315 degrees from 1 -1 center at 1 0
\endpicture
%
\break
{\small\ifx\thesamplenumber\relax 
Volume \else Sample
\fi\thevolumenumber\ (\thevolumeyear)
\startpage--\finishpage\nl
Published: \publishdate}
\vglue 0.5truein plus 0.4fil minus 0.1truein

{\parskip=0pt\leftskip 0pt plus 1fil\def\\{\par\smallskip}{\ifplaintex\large
\else\Large\fi\bf\thetitle}\par\medskip}   

\vglue 0pt plus 0.1fil 

{\parskip=0pt\leftskip 0pt plus 1fil\def\\{\par}{\sc\theauthors}
\par\medskip}

\vglue 0pt plus 0.1fil 

{\small\parskip=0pt\let\newline\\
{\leftskip 0pt plus 1fil\def\\{\par}{\sl\theaddress}\par}
\expandafter\ifx\theemail\relax    
\relax\else\vglue 5pt plus 0.02fil minus 2pt\def\\{\stdspace{\rm 
and}\stdspace} 
\cl{Email:\stdspace\tt\theemail}\fi
\ifx\theurl\relax                  
\relax\else\vglue 5pt plus 0.02fil minus 2pt\def\\{\stdspace{\rm 
and}\stdspace}
\cl{URL:\stdspace\tt\theurl}\fi\par}

\vglue 7pt plus 0.3fil minus 3pt

{\bf Abstract}
\vglue 5pt plus 0.1fil minus 2pt

\theabstract

\vglue 7pt plus 0.3fil minus 3pt

{\bf AMS Classification numbers}\quad Primary:\quad \theprimaryclass

Secondary:\quad \thesecondaryclass

\vglue 5pt plus 0.3fil minus 2pt

{\bf Keywords}\quad \thekeywords

\vglue 10pt plus 0.5fil minus 5pt

{\small  Proposed: \theproposer\hfill Received: \receiveddate\nl
Seconded: \theseconders\hfill 
\ifx\reviseddate\relax                         
Accepted: \accepteddate                        
\else
Revised: \reviseddate                          
\fi}
\eject
}       

\let\maketitlepage\maketitlep
\let\maketitle\maketitlepage


\font\phead=cmsl9 scaled 950
\font\lhead=cmsl9 scaled 1050
\font\pnum=cmbx10 scaled 913
\font\lnum=cmbx10 
\font\pfoot=cmsl9 scaled 950
\font\lfoot=cmsl9 scaled 1050
\ifplaintex
\headline{\vbox to 0pt{\vskip -4.5mm\line{\small\phead\ifnum
\count0=\startpage ISSN 1364-0380 (on line)
1465-3060 (printed) \hfill {\pnum\folio}\else\ifodd\count0\def\\{ }%
\ifx\theshorttitle\relax\thetitle\else\theshorttitle\fi\hfill{\pnum\folio}
\else\def\\{ and }{\pnum\folio}\hfill\ifx\theshortauthors\relax\theauthors
\else\theshortauthors\fi\fi\fi}\vss}}
\footline{\vbox to 0pt{\vglue 0mm\line{\small\pfoot\ifnum\count0=\startpage
\copyright\ \gtp\hfill\else
\gt, Volume \thevolumenumber\ (\thevolumeyear)\hfill\fi}\vss
}}
\else
\makeatletter
\def\@oddhead{{\small\lhead\ifnum\count0=\startpage ISSN 1364-0380 (on line)
1465-3060 (printed) \hfill {\lnum\number\count0}\else\ifodd\count0
\def\\{ }\ifx\theshorttitle\relax \thetitle \else\theshorttitle\fi\hfill
{\lnum\number\count0}\else\def\\{ and }{\lnum\number\count0}
\hfill\ifx\theshortauthors\relax 
\theauthors\else\theshortauthors\fi\fi\fi}}\def\@evenhead{\@oddhead}
\def\@oddfoot{\small\lfoot\ifnum\count0=\startpage\copyright\ \gtp\hfill\else
\gt, Volume \thevolumenumber\ (\thevolumeyear)\hfill\fi}
\def\@evenfoot{\@oddfoot}
\makeatother
\fi


\newwrite\gtoutfile
\long\gdef\makeheadfile{  
{\def\\{, }\def\s{ }
\immediate\openout\gtoutfile head.xxx
\immediate\write\gtoutfile{To: math@arxiv.org}
\immediate\write\gtoutfile{Subject: put or rep NNNNN:pppp}
\immediate\write\gtoutfile{--text follows this line--}
\immediate\write\gtoutfile{Proxy-for: \ifx\theasciiauthors\relax
\theauthors\else\theasciiauthors\fi\s<\ifx\theasciiemail\relax\theemail\else\theasciiemail\fi>}
\immediate\write\gtoutfile{\noexpand\\}
\immediate\write\gtoutfile{Authors: \ifx\theasciiauthors\relax
\theauthors\else\theasciiauthors\fi}
\immediate\write\gtoutfile{Title: \ifx\theasciititle\relax
\thetitle\else\theasciititle\fi}
\immediate\write\gtoutfile{Subj-class: GT or SG or MG etc}
\immediate\write\gtoutfile{MSC-class: \theprimaryclass\ifx\thesecondaryclass\relax\else, \thesecondaryclass\fi}
\immediate\write\gtoutfile{Journal-ref: Geom. Topol. \thevolumenumber
(\thevolumeyear) \startpage-\finishpage}
\immediate\write\gtoutfile{Comments: Published by Geometry and Topology at}
\immediate\write\gtoutfile{\s\s http://www.maths.warwick.ac.uk/gt/GTVol\thevolumenumber/paper\thepapernumber.abs.html}
\immediate\write\gtoutfile{\noexpand\\}
\immediate\write\gtoutfile{}
\ifx\theasciiabstract\relax
\immediate\write\gtoutfile{\theabstract}\else
\immediate\write\gtoutfile{\theasciiabstract}\fi
\immediate\write\gtoutfile{}
\immediate\write\gtoutfile{\noexpand\\}
\immediate\write\gtoutfile{}
\immediate\closeout\gtoutfile}}  

\def\maketitlepage{\maketitlep\makeheadfile}
\let\maketitle\maketitlepage

\lognumber{150}
\volumenumber{5}\papernumber{17}\volumeyear{2001}
\pagenumbers{521}{550}
\revised{11 May 2001}
\accepted{18 May 2001}
\proposed{Walter Neumann}
\seconded{Jean-Pierre Otal, Steve Ferry}
\received{23 November 2000}
\published{21 May 2001}

\let\cal\Cal      
\let\\\par
\def\topmatter{\relax}
\def\endtopmatter{\maketitlepage}
\let\gttitle\title
\def\title#1\endtitle{\gttitle{#1}}
\let\gtauthor\author
\def\author#1\endauthor{\gtauthor{#1}}
\let\gtaddress\address
\def\address#1\endaddress{\gtaddress{#1}}
\def\affil#1\endaffil{\gtaddress{#1}}
\let\gtemail\email
\def\email#1\endemail{\gtemail{#1}}
\def\subjclass#1\endsubjclass{\primaryclass{#1}}
\let\gtkeywords\keywords
\def\keywords#1\endkeywords{\gtkeywords{#1}}
\def\heading#1\endheading{{\def\S##1{\relax}\def\\{\relax\ignorespaces}
    \section{#1}}}
\def\head#1\endhead{\heading#1\endheading}

\def\subhead#1\endsubhead{\sh{#1}}
\def\subsubhead#1\endsubsubhead{\sh{#1}}
\def\specialhead#1\endspecialhead{\sh{#1}}
\def\demo#1{\rk{#1}\ignorespaces}
\def\enddemo{\ppar}
\let\remark\demo
\def\endremark{}
\let\definition\demo
\def\enddefinition{\ppar}

\def\qed{\ifmmode\quad\sq\else\hbox{}\hfill$\sq$\par\goodbreak\rm\fi}  
\def\proclaim#1{\rk{#1}\sl\ignorespaces}
\def\endproclaim{\rm\ppar}
\def\cite#1{[#1]}
\newcount\itemnumber
\def\roster{\items\itemnumber=1}
\def\endroster{\enditems}
\let\itemold\item
\def\item{\itemold{{\rm(\number\itemnumber)}}%
\global\advance\itemnumber by 1\ignorespaces}
\def\S{section~\ignorespaces}  
\def\date#1\enddate{\relax}
\def\thanks#1\endthanks{\relax}   
\def\dedicatory#1\enddedicatory{\relax}  
\let\footnote\plainfootnote

\def\Refs{\ppar{\large\bf References}\ppar\bgroup\leftskip=25pt
\frenchspacing\parskip=3pt plus2pt\small}       
\def\endRefs{\egroup}
\def\widestnumber#1#2{\relax}
\def\endrefitem{}
\def\refdef#1#2#3{\def#1{\leavevmode\unskip\endrefitem#2\def\endrefitem{#3}}}
\def\ref{\par}
\def\endref{\endrefitem\par\def\endrefitem{}}
\refdef\key{\noindent\llap\bgroup[}{]\ \ \egroup}
\refdef\no{\noindent\llap\bgroup[}{]\ \ \egroup}
\refdef\by{\bf}{\rm, }
\refdef\manyby{\bf}{\rm, }
\refdef\paper{\it}{\rm, }
\refdef\book{\it}{\rm, }
\refdef\jour{}{ }
\refdef\vol{}{ }
\refdef\yr{(}{) }
\refdef\ed{(}{, Editor) }
\refdef\publ{}{ }
\refdef\inbook{from: ``}{'', }
\refdef\pages{}{ }
\refdef\page{}{ }
\refdef\paperinfo{}{ }
\refdef\bookinfo{}{ }
\refdef\publaddr{}{ }
\refdef\moreref{}{ }
\refdef\eds{(}{, Editors)}
\refdef\bysame{\hbox to 3 em{\hrulefill}\thinspace,}{ }
\refdef\toappear{(to appear)}{ }
\refdef\issue{no.\ }{ }
\newcount\refnumber\refnumber=1
\def\refkey#1{\expandafter\xdef\csname cite#1\endcsname{\number\refnumber}%
\global\advance\refnumber by 1}
\def\cite#1{[\csname cite#1\endcsname]}
\def\Cite#1{\csname cite#1\endcsname}  
\def\key#1{\noindent\llap{[\csname cite#1\endcsname]\ \ }}
\refkey{Ba}
\refkey{BB1}
\refkey{BB2}
\refkey{BH}
\refkey{Bo}
\refkey{Br}
\refkey{CD}
\refkey{D}
\refkey{G}
\refkey{K}
\refkey{KL}
\refkey{N}
\refkey{R}
\refkey{T}

\define\s{\Bbb S}
\define\r{\Bbb R}
\define\g{\gamma}
\define\a{\Cal A}
\define\lk{\text{lk}\,}
\define\sig{\Sigma_W}
\define\st{\text{st}\,(x)}

\topmatter
\title Metric characterizations of spherical\\and Euclidean buildings
\endtitle
\author Ruth Charney\\Alexander Lytchak \endauthor

\address
Mathematics Department, Ohio State University\\231 
W 18th Ave, Columbus, OH 43210, USA\\\smallskip {\rm and}\\\smallskip
\\Mathematisches 
Institut der Universit\"at Bonn\\Wegelerstra\ss e 10, D-53115 Bonn, Germany
\endaddress
\asciiaddress{
Mathematics Department, Ohio State University\\231 
W 18th Ave, Columbus, OH 43210, USA\\Mathematisches 
Institut der Universitat Bonn\\Wegelerstrasse 10, D-53115 Bonn, Germany}

\gtemail{charney@math.ohio-state.edu\qua {\rm and}\qua 
lytchak@math.uni-bonn.de}
\asciiemail{charney@math.ohio-state.edu, lytchak@math.uni-bonn.de}

\abstract
A building is a simplicial complex with a covering by Coxeter
complexes (called apartments) satisfying certain combinatorial
conditions. A building whose apartments are spherical (respectively 
Euclidean) Coxeter complexes has a natural piecewise spherical 
(respectively Euclidean) metric with nice geometric properties. We show 
that spherical and Euclidean buildings are completely characterized 
by some simple, geometric properties.
\endabstract

\primaryclass{20E42}

\secondaryclass{20F65}

\gtkeywords{Buildings, CAT(0) spaces, spherical buildings, Euclidean buildings,
metric characterisation}

\endtopmatter

\document
\sectionnumber-1

\head Introduction \endhead

In recent years, much attention has been given to curvature properties
of piecewise Euclidean and piecewise spherical complexes.  A notion of
curvature bounded above for such complexes was introduced by
Alexandrov in the 1950's and further developed in the 1980's by
Gromov. The curvature bound is defined as a condition on the shape of
triangles (they must be sufficiently ``thin'') and is known as a
CAT--inequality (Comparison inequality of Alexandrov--Toponogov). Some
particularly nice examples of spaces satisfying CAT--inequal\-it\-ies are
spherical and Euclidean buildings which come equipped with a natural
piecewise spherical or Euclidean metric.

Buildings also satisfy other nice metric properties. A spherical
building $X$, for example, is easily seen to have diameter $\pi$, as
does the link of any simplex in $X$. It is natural to ask whether
Euclidean and spherical buildings are characterized by their metric
properties. In this paper, we give several metric characterizations of
buildings. For example we prove

\proclaim{Theorem} Let $X$ be a connected, piecewise spherical 
(respectively Euclidean) complex of dimension $n \geq 2$ satisfying
\roster
\item $X$ is CAT(1) (respectively CAT(0)).
\item Every $(n-1)$--cell is contained in at least two $n$--cells.
\item Links of dimension $\geq 1$ are connected.
\item Links of dimension $1$ have diameter $\pi$.
\endroster
Then $X$ is isometric to a spherical building (respectively a metric
Euclidean building).
\endproclaim

Metric Euclidean buildings are products of irreducible Euclidean
buildings, cones on spherical buildings, trees, and nonsingular
Euclidean spaces (see Section 5).

Another metric property of buildings is that for every local geodesic
$\g$, the set of directions in which $\g$ can be geodesically
continued is non-empty and discrete. We call this the ``discrete
extension property''.

\proclaim{Theorem}  Let $X$ be a connected, piecewise spherical 
(respectively Euclidean) complex of dimension $n \geq 2$ satisfying
\roster
\item $X$ is CAT(1) (respectively CAT(0)).
\item $X$ has the discrete extension property.
\endroster
Then $X$ is isometric to a spherical building (respectively metric Euclidean
building).
\endproclaim

Werner Ballmann and Michael Brin, studying the question of rank
rigidity for piecewise Euclidean complexes of nonpositive curvature,
have obtained related results in \cite{BB1} and \cite{BB2}, including
a metric characterization of spherical and Euclidean buildings of
dimension 2.  Bruce Kleiner has also described (unpublished) a metric
characterization of Euclidean buildings under the assumption that
every geodesic is contained in an $n$--flat.

The first two sections of the paper contain background about buildings
and geodesic metric spaces. The key problem in identifying a building
is the construction of enough apartments. Sections 3, 4, and 5 are
devoted to this task. Section 6 considers the 1--dimensional case, and
Section 7 combines these results to arrive at the main theorems.

The first author would like to thank Bruce Kleiner for helpful
conversations. The second author is indebted to Werner Ballmann for
suggesting the question and for his continued support during the
development of this paper.  The first author was partially supported
by NSF grant DMS-9803374.

\head Buildings \endhead

In this section we review some definitions and terminology. For more
details about buildings, see \cite{Br} and \cite{Bo}.

Let $S$ be a finite set. A {\it Coxeter matrix} on $S$ is a symmetric
function $m \co  S \times S \to \{1,2,\dots ,\infty \} $ such that
$m(s,s)=1$ and $m(s,t) \geq 2$ for $s \neq t$. The {\it Coxeter group}
associated to $m$ is the group $W$ given by the presentation
$$ W=\langle ~S ~\vert ~(st)^{m(s,t)} = 1,\, s,t \in S ~\rangle. $$
The pair $(W,S)$ is called a Coxeter system.
If $T \subset S$, then the subgroup $W_T$ generated by $T$ is the
Coxeter group associated to $m \vert_{T \times T}$. The Coxeter system
$(W,S)$ is {\it irreducible} if there is no non-trivial partition $S=S_1
\amalg S_2$ such that $W=W_{S_1} \times W_{S_2}$.

The Coxeter group $W$ can be realized as a discrete group of linear
transformations of an $n$--dimensional vector space $V$, with the
generators $s \in S$ acting as reflections across the walls of a
simplicial cone. This action preserves a bilinear form $B$ on $V$
represented by the matrix $B(s,t)= -\text{cos}\,(\frac{\pi}{m(s,t)})$
(where $\frac{\pi}{\infty}$ is taken to be 0). $W$ is finite if and only if
this form is positive definite. In this case, $W$ acts as a group of
orthogonal transformations and the action restricts to the unit sphere
$S(V)$ in $V$. Hence $W$ is called a {\it spherical Coxeter group}.

If $B$ is positive semi-definite (but not definite) and $W$ is
irreducible, then the action of $W$ on $V$ induces an action on an
(n-1)-dimensional affine space $\Bbb R^{n-1} = V / V^\perp$ with the
generators acting as affine reflections across the walls of a simplex.
In this case, $W$ is called an {\it irreducible Euclidean Coxeter
group.} A key fact about irreducible Euclidean Coxeter groups is that
for any proper subgroup $T \subset S$, $W_T$ is a spherical Coxeter
group.  More generally, we call $W$ a {\it Euclidean Coxeter group} if
it is a direct product of irreducible Euclidean Coxeter groups.

To any Coxeter group $W$, one can associate a simplicial complex
$\sig$ called the {\it Coxeter complex} for $W$. In the case of
spherical and Euclidean Coxeter groups, the Coxeter complex has a
simple, geometric description. Let $M=\s^{n-1}=S(V)$ if $W$ is
spherical or $M=\Bbb R^{n-1}=V/V^{\perp}$ if $W$ is irreducible
Euclidean. For each element $r \in W$ which acts as a reflection on
$M$ (namely, $r$ is a generator or conjugate of a generator), $r$
fixes some hyperplane, called a {\it wall} of $M$.  The walls divide
$M$ into simplices. The resulting simplicial complex, $\sig$, is the
Coxeter complex for $W$. If $W$ is a product of irreducible Euclidean
Coxeter groups, then $\sig$ is the product of the corresponding Coxeter
complexes. The top dimensional simplices (or cells) of $\sig$ are called {\it
chambers}. $W$ acts freely transitively on the set of chambers of
$\sig$ and the stabilizer of any lower dimensional cell $\sigma$ is
conjugate to $W_T$ for some $T \subset S$.

There are several equivalent definitions of buildings. The most
convenient for our purposes is the following (see \cite{Br}).

\definition{Definition 1.1} A {\it building} is a simplicial complex
$X$ together with a collection of subcomplexes $\a$, called
{\it apartments}, satisfying
\roster
\item each apartment is isomorphic to a Coxeter complex,
\item any two simplices of $X$ are contained in a common apartment,
\item if two apartments $A_1, A_2$ share a chamber, then there is an
isomorphism $A_1 \to A_2$ which fixes $A_1 \cap A_2$ pointwise.
\endroster
If, in addition, every codimension 1 simplex is contained in at least
three chambers, then $X$ is a {\it thick} building.
It follows from conditions (2) and (3) that all of the apartments are
isomorphic to the same Coxeter complex $\sig$. We say that a building
$X$ is {\it spherical} (respectively {\it Euclidean}) if $W$ is spherical
(respectively Euclidean).
\enddefinition

Although the collection of apartments $\a$ is not, in general
unique, there is a unique maximal set of apartments. We will always
assume $\a$ to be maximal.

If $X_1, X_2$ are spherical buildings with associated Coxeter groups
$W_1,W_2$, then the join $X_1 * X_2$ is a spherical building with
Coxeter group $W_1 \times W_2$. In particular, the suspension $\Sigma
X_2 = \s^0 * X_2$ is a building with Coxeter group $\Bbb Z/2 \times
W_2$. (Note that for $k>0$, the simplicial structure on $\s^k * X_2$
depends on a choice of identification of $\s^k$ with the (k+1)-fold
join $\s^0* \cdots *\s^0$. Despite this slight ambiguity, we will
consider the $k$--fold suspension of a building to be a building.)
Conversely, if $X$ is a spherical building whose Coxeter group
$W$ splits as a product $W_1 \times W_2$, then $X$ can be decomposed as
the join of a building for $W_1$ and a building for $W_2$ (see
\cite{R}, Theorem 3.10). Similarly, any Euclidean building splits as a
product of irreducible Euclidean buildings.

\head Metrics \endhead

A metric space $(X,d)$ is a {\it geodesic metric space} if for any two
points $x,y \in X$, there is an isometric embedding $\g \co  [0,a] \to X$
with $\g(0)=x$ and $\g(a)=y$. Such a path is called a {\it geodesic
segment} or simply a {\it geodesic}
from $x$ to $y$. An isometric embedding of $\Bbb R$ into $X$ is also
called a {\it geodesic}, and an isometric embedding of $[0,\infty)$ is
called a {\it ray.}

A {\it piecewise Euclidean} (respectively {\it piecewise spherical}) complex
is a polyhedral cell complex $X$ together with a metric $d$ such that
each cell of $X$ is isometric to a convex polyhedral cell in $\r^n$
(respectively $\s^n$) for some $n$, and 
$$ d(x,y) = \text{inf \{ length} (\g) \mid \g \text{ is a path from
$x$ to $y$} \}$$ for any $x,y \in X$. We will also assume that the
metric $d$ is a complete, geodesic metric. In particular, the infimum
$d(x,y)$ is realized by the length of some path $\g$.

If $X$ is a piecewise spherical or Euclidean complex and $x$ is a point
in $X$, then the set of unit tangent vectors to $X$ at $x$ is called the
{\it link} of $x$ and denoted $\lk (x,X)$. It comes equipped with the
structure of a piecewise spherical complex, since the link of $x$ in
a single $n$--cell is isometric to a polyhedral cell in $\s^{n-1}$. If
$\sigma$ is a $k$--cell in $X$, we define $\lk (\sigma, X)$ to be the
set of unit tangent vectors orthogonal to $\sigma$ at any point $x$
in the relative interior of $\sigma$. This set also has a natural
piecewise spherical structure and we can identify
$$\lk (x,X) = \lk (x,\sigma) * \lk (\sigma , X) = \s^{k-1} * \lk
(\sigma ,X) $$ where the joins $*$ are orthogonal joins in the sense
of \cite{CD}.  (See \cite{BH} or the appendix of \cite{CD} for a
discussion of joins of piecewise spherical complexes.)

In some cases, we may wish to consider spaces which do not have a
globally defined cell structure. For this, we introduce the notion of
a locally spherical space of dimension $n$. The definition is
inductive. A locally spherical space of dimension 0 is a nonempty
disjoint union of points. A {\it locally spherical (respectively Euclidean)
space of dimension $n$}, $n>0$, is a complete geodesic metric space
$(X,d)$ for which every point $x$ has a neighborhood isometric to a
spherical (respectively Euclidean) cone on a locally spherical space $L_x$ of
dimension $n-1$. We call such a neighborhood a {\it conelike
neighborhood} of $x$. Clearly, $L_x= \lk (x,X)$. A piecewise spherical
(respectively Euclidean) complex is a locally spherical (respectively Euclidean)
space of dimension $n$ if and only if every cell is contained in an
$n$--dimensional cell.

The basis for our metric characterization of buildings will be the
CAT--inequal\-it\-ies defined by Gromov in \cite{G}.
Let $(X,d)$ be a complete, geodesic metric space and
let $T$ be a geodesic triangle in $X$. A Euclidean comparison
triangle for $T$ is a triangle $T'$ in $\r^2$ with the same side
lengths as $T$. We say $X$ is a {\it CAT(0) space} if every geodesic
triangle $T$ is ``thin" relative to its comparison triangle $T'$. That
is, given any points $x,y \in T$, the distance from $x$ to $y$ in $X$
is less than or equal to the distance in $\r^2$ between the
corresponding points $x',y' \in T'$. We define a {\it CAT(1) space}
similarly by comparing geodesic triangles $T$ in $X$ with spherical
comparison triangles $T'$ in $\s^2$. In this case, however, we only
require the thinness condition to hold for triangles $T$ of perimeter
$\leq \pi$ (since no comparison triangle exists with perimeter $>
\pi$). 

In the next two theorems we collect some facts about CAT(0) and CAT(1)
spaces. These are due to Gromov, Ballmann, Bridson and others. A good
source of proofs is \cite{BH} or \cite{Ba}.

\proclaim{Theorem 2.1} Let $X$ be a piecewise (or locally) Euclidean
geodesic metric space.
\roster
\item $X$ is locally CAT(0) if and only if $\lk (\sigma ,X)$ is CAT(1)
for every cell $\sigma$.
\item $X$ is CAT(0) if and only if it is locally CAT(0) and simply
connected. 
\item If $X$ is CAT(0), then any two points in $X$ are connected by a
unique geodesic and any path which is locally geodesic, is geodesic.
\endroster\endproclaim

\proclaim{Theorem 2.2} Let $X$ be a piecewise (or locally) spherical
geodesic metric space.
\roster
\item $X$ is locally CAT(1) if and only if $\lk (\sigma ,X)$ is CAT(1)
for every cell $\sigma$.
\item $X$ is CAT(1) if and only if it is locally CAT(1), any two
points of distance $< \pi$ are connected by a unique geodesic, and
these geodesics vary continuously with their endpoints.
\item If $X$ is CAT(1), then  any path of length $\leq \pi$ 
which is locally geodesic, is geodesic.
\endroster\endproclaim

A Euclidean (respectively spherical) building $X$ of dimension $n$ comes
equipped with a natural piecewise Euclidean (respectively piecewise
spherical) metric in which each apartment is isometric to $\Bbb R^n$
(respectively $\s^n$) with the Coxeter group $W$ acting by
isometries. In the spherical case, there is a unique such metric. In
the Euclidean case, this metric is determined only up to scalar
multiple on each irreducible factor. {\it In this paper,
for a piecewise Euclidean (respectively spherical) complex $(X,d)$, the
statement that $X$ is a Euclidean (respectively spherical) building will mean
that, the cell structure on $X$ satisfies the conditions of Definitions 1.1
and that the metric on $X$ is the natural building metric.}

Define the diameter of $X$ to be
$$diam(X)= \text{sup\,}\{d(x,y) \vert x,y \in X \}. $$
The natural metric on a building satisfies a number of nice properties
which are described in the proposition below. 

\proclaim{Proposition 2.3} Let $X$ be a Euclidean (respectively spherical)
building of dimension $n$ with the natural metric. Then:
\roster
\item $X$ is CAT(0) (respectively $X$ is CAT(1) and
$\text{daim\,}(X)=\pi$). 
\item For any simplex $\sigma$ of codimension $\geq 2$, $\lk
(\sigma,X)$ is a spherical building. In particular, $\lk (\sigma,X)$ 
is CAT(1) and $\text{diam\,}(\lk (\sigma,X)) = \pi.$
\item A subspace $A \subset X$ is an apartment if and only if the
intrinsic metric on $A$ is isometric to $\r^n$ (respectively $\s^n$).
Moreover, the inclusion $A \hookrightarrow X$ is an isometric embedding.
\endroster\endproclaim

\demo{Proof} (1) and (2) are well known. (1) follows from the fact
that every geodesic in $X$ is contained in an apartment (see \cite{D} or
\cite{Br}). (2) follows from the fact that the isotropy group of a 
simplex $\sigma$ in
$\sig$ is a spherical Coxeter group $W_{\sigma}$. The action of
$W_{\sigma}$ on the sphere $\lk (\sigma, \sig)$ gives a natural
identification of $\lk (\sigma, \sig)$ with the Coxeter complex for
$W_\sigma$. These constitute the apartments of $\lk (\sigma, X)$.

For (3), let $M=\r^n$ if $X$ is Euclidean and $M=\s^n$ if $X$ is
spherical. Consider the collection of subspaces
$$ \a = \{ A \subset X \vert \text{$A$ is isometric to $M$} \}. $$ 
By definition of the metric on $X$, $\a$ contains all the apartments
of $X$. Since we are assuming the set of apartments to be maximal, it
suffices to show that $\a$ satisfies the conditions in Definition 1.1
for a system of apartments.

Observe first that any subspace $A$ isometric to $M$ is necessarily a
subcomplex of $X$ since its intersection with any $n$--simplex
$\sigma$ must be both open and closed in $\sigma$. 
By induction on the dimension of $X$, we may assume that for any
simplex $\sigma \subset A$, $\lk (\sigma ,A) \hookrightarrow \lk
(\sigma ,X)$ is an isometric embedding. It follows that the embedding
$A \hookrightarrow X$ preserves local geodesics. If $X$ is Euclidean,
then it is CAT(0), hence local geodesics are geodesics. If $X$ is
spherical, then it is CAT(1), hence local geodesics of length $\leq
\pi$ are geodesics. In either case, we conclude that $A
\hookrightarrow X$ maps geodesics to geodesics, so it is an isometric
embedding. 

Let $A \in \a$ and let $\sigma \subset A$ be an $n$--simplex. Fix an
isometry $\alpha_0$ of $\sigma$ with the fundamental chamber
$\sigma_0$ of a Coxeter complex $\sig$. This isometry extends uniquely
to an isometry $\alpha \co  A \to \sig$. Since every $n$--simplex in $A$
is isometric to $\sigma_0$, $\alpha$ is also a simplicial
isomorphism. Thus, $A$ is a Coxeter complex. Moreover, if $A' \in \a$
also contains $\sigma$, and $\alpha ' \co  A' \to \sig$ is an isometry
extending $\alpha_0$, then $\alpha^{-1} \circ \alpha ' \co  A \to A'$ is
an isometry fixing $\sigma$ and hence fixing all of $A \cap A'$.
Finally, since $\a$ contains a system of apartments, any two simplices
of $X$ are contained in some $A \in \a$. Thus, $\a$ satisfies the
conditions for a system of apartments.
\qed\enddemo

\head Spherical buildings \endhead

In this section we prove a partial converse to Proposition 2.3. It
will form the inductive step to one of the main theorems in Section
7.
 
\proclaim{Theorem 3.1} Suppose $X$ is a connected, piecewise
spherical cell complex of dimension $n \geq 2$ satisfying
\roster
\item $X$ is CAT(1),
\item $\lk (x,X)$ is a spherical building for every vertex $x \in X$.
\endroster
Then $X$ is a spherical building.
\endproclaim

Before embarking on the proof, we make several observations about the
hypotheses. First, If $\sigma$ is any cell in $X$, then $L_\sigma =
\lk (\sigma,X)$ is a spherical building. For if $v$ is any vertex of 
$\sigma$, then
$$L_\sigma =\lk (\lk(v,\sigma),\lk (v,X))$$ 
and $\lk (v,X)$ is a spherical building by hypothesis. Thus, it
follows from Proposition 2.3(2) that $L_\sigma$ is a spherical
building. Moreover, if $x$ is any point in $X$, not necessarily a
vertex, then $\lk (x,X)$ is also a spherical building. For if $x$ lies
in the relative interior of a $k$--cell $\sigma$, then $\lk (x,X) =
\s^{k-1} * L_\sigma = \Sigma^k (L_\sigma)$.

By Proposition 2.3, there is an obvious candidate for a system of
apartments for $X$, namely
$$ \a = \{ A \subset X \vert \text{$A$ is isometric to $\s ^n$} \}.$$
As in the proof of Proposition 2.3(3), it is easy to show that any
such subspace $A$ is a subcomplex of $X$ and the
inclusion $A\hookrightarrow X$ is an isometric embedding.

The key problem in the proof of Theorem 3.1 is to construct
enough of these subcomplexes. The idea is as follows. For any pair of
antipodal points (two points are {\it antipodal} if they have distance
$\pi$) and any apartment $A_x$ in $\lk (x,X)$, we construct an
apartment $A$ in $X$ by propagating geodesics from $x$ to $y$ in every
direction in $A_x$.

We begin with a key technical lemma.
Some additional notation will be needed for the proof. If $x \in X$
and $\g$ is a geodesic emanating from $x$, let $\g_x \in \lk (x,X)$
denote the tangent vector to $\g$ at $x$. Let $\st$ denote the closed
star of $x$, that is, $\st$ is the union of the closed simplices
containing $x$. (In the locally spherical context, $\st$ will denote
a conelike polyhedral neighborhood of $x$.)

For $\sigma$ a spherical $(n-1)$--cell, the
spherical suspension $\Sigma (\sigma)$, viewed as a subspace of $\s^n$, is
called a {\it spherical sector}. When $n=2$, it is also called a {\it
spherical lune}.

We prove the next lemma under slightly
more general hypotheses for use in the next section. In particular, we
do not assume that $X$ is globally CAT(1).

\proclaim{Lemma 3.2} Suppose $X$ is a locally spherical space of
dimension $n \geq 2$ such that the link of every point in $X$ is
isometric to a building. Let $\g$ be a local geodesic of length $\pi$
from $x$ to $y$ and let $A_x$ be an apartment in $L_x$ containing
$\g_x$. Then there is a neighborhood $N_x \subset A_x$ of $\g_x$ and a
unique locally isometric map $F$ of the spherical sector 
$\Sigma (N_x)$ into $X$ such that
\roster
\item for any $v \in N_x$, the restriction of $F$ to $\Sigma (v)$
$(=\s^0 * \{ v\})$ is a
local geodesic from $x$ to $y$ with tangent vector $v$, and
\item the restriction of $F$ to $\Sigma (\g_x)$ is precisely $\g$.
\endroster
\endproclaim

\demo{Proof} Divide $\g$ into segments $\g_1,\g_2, \dots ,\g_k$ with
endpoints $x=x_0, x_1, \dots ,x_k$ such that each $\g_i$
lies in $\text{st}\,(x_{i-1})$.
Let $N_x$ be an $\epsilon$--ball in $A_x$ centered
at $\g_x$ and let $S$ be the spherical sector $S=\Sigma (N_x)$.
 
For each vector $v \in N_x$, there is a unique geodesic
segment $\g_1^v$ in $\st$ from $x$ to $\partial(\st )$ whose tangent
at $x$ is $v$.
Let $B^1$ be the subspace of $\st$ consisting of the union of these 
geodesic segments. Identifying $\g_1^v$ with an initial segment of
$\Sigma (v)$ gives an isometry $F_1$ of a polyhedral
subspace of  $S$  onto $B_1$.

Next, consider the $(n-1)$--dimensional building $L_1 = \lk(x_1, X)$.
The tangent vectors to $\g_1$ and $\g_2$ at $x_1$ form a pair of
antipodal points $a_1,a_2$ in $L_1$ (since the concatenation $\g_1
\cdot \g_2$ is geodesic), and $\lk (x_1, B_1)$ is a neighborhood of
$a_1$ in $L_1$ isometric to a spherical $(n-1)$--cell. The union of
geodesics in $L_1$ from $a_1$ to $a_2$ with an initial segment lying
in this spherical cell forms an apartment $A_1 \subset L_1$. The
geodesic segments emanating from $x_1$ in directions $A_1$ form a
spherical $n$--cell $C$ in $\text{st\,}(x_1)$. Shrinking the original
$\epsilon$--neighborhood $N_x$ if necessary, we may assume that all of
the segments $\g_1^v$ end in $C$. There is then a unique locally
geodesic continuation of $\g_1^v$ across $C$, ending in $\partial
(\text{st\,}(x_1))$. Call this new segment $\g_2^v$. Let $B_2$ be the
union of the local geodesics $\g_1^v
\cdot \g_2^v, \, v \in N_x$. Then $F_1$ extends in an obvious
manner to a local isometry $F_2$ from a polyhedral subspace of 
$S$ onto $B_2$.

We repeat this process at each $x_i$ until we get geodesics $\g^v =
\g_1^v \cdot \g_2^v \cdots \g_k^v$ for every $v \in N_x$ and a local 
isometry $F$ from $S$ onto $B_k = \bigcup \g^v$ as required.
\qed\enddemo

Returning to the hypothesis of Theorem 3.1, we can now construct
apartments in $X$.

\proclaim{Lemma 3.3} Let $X$ be as in Theorem 3.1. Suppose $x,y \in X$ 
are antipodal points and $A_x$ is an apartment in $L_x = \lk (x,X)$. 
Then the following hold.
\roster
\item For every $v \in A_x$, there exists a unique geodesic $\g^v$
from $x$ to $y$ whose tangent vector at $x$ is $v$.
\item The union of all such $\g^v$, $v \in A_x$, is isometric to 
$\s^n$.
\endroster\endproclaim

\demo{Proof} First note that since $X$ is CAT(1), any local geodesic
of length $\leq \pi$ is a geodesic. Moreover, if two geodesics from 
$x$ to $y$, $\g$
and $\g '$, have the same tangent vectors $\g_x =\g'_x =v$, then they
must agree inside $\st$, hence they must agree everywhere.
(Otherwise, we get a geodesic digon of length $<2\pi$.) Thus,
geodesics from $x$ to $y$ are uniquely determined by their tangents
at $x$.

Consider the set $C=\{ v\in A_x \vert \text{$\g^v$ exists}\}$. By
Lemma 3.2, this set is open in $A_x$. We claim that it is also closed.
To see this, first note that if $v_1,v_2 \in C$ are
points of distance $\alpha$ in $A_x$, then for any $t \in [0,\pi]$,
$d(\g^{v_1}(t), \g^{v_2}(t)) \leq \alpha$. This can be seen by
comparing the digon formed by $\g^{v_1}, \g^{v_2}$ in $X$ with a
digon $\g'_1, \g'_2$ of angle $\alpha$ connecting a pair of antipodal
points $x',y'$ in $\s^2$. Inside $\st$, these two digons are
isometric. That is, for sufficiently small $\epsilon$, the distance
from $z_1=\g^{v_1}(\epsilon)$ to $z_2=\g^{v_2}(\epsilon)$ in $X$ is
the same as the distance between the corresponding points $z'_1,z'_2$
in $\s^2$. Thus, $z'_1,z'_2,y'$ is a spherical comparison triangle
for $z_1,z_2,y$. It follows from the CAT(1) condition that the
distance between $\g^{v_1}(t)$ and $\g^{v_2}(t)$ in $X$ is less than
or equal to the corresponding distance in $\s^2$ for all $t$. 
In particular, if
$(v_i)$ is a sequence of points in $A_x$ converging to $v$ with
$v_i \in C$, then $(\g^{v_i})$ converges uniformly to a
path $\g$ from $x$ to $y$ with $\g_x=v$. This path has length $\pi$
since each $\g^{v_i}$ has length $\pi$. Hence $\g$ is geodesic and $v
\in C$.

Since $C$ is both open and closed, it is either empty or all of
$A_x$. Since $X$ is a geodesic metric space, there must exist at least
one geodesic $\g$ from $x$ to $y$. If $v$ is any point in $L_x$, then
there exists an apartment $A_x$ containing both $\g_x$ and $v$. For
this apartment, $C$ is nonempty, hence $v \in C$. Thus, there is a
geodesic $\g^v$ with tangent vector $v$ as desired. This proves (1).

By Lemma 3.2, we know that the map $F\co  \Sigma (A_x) \to X$ taking
$\Sigma (v)$ to $\g^v$ is locally isometric. Since local geodesics
of length $\leq \pi$ are geodesic in $X$, this map is an isometry onto
its image. This proves (2).
\qed\enddemo

The spheres constructed in Lemma 3.3 give us a large number of
apartments. It is now easy to show that $X$ is a building.

\proclaim{Lemma 3.4} $X$ has diameter $\pi$.
\endproclaim

\demo{Proof} Since every point in $\lk (x,X)$ has an antipodal point,
geodesics in $X$ are locally extendible. Since $X$ is CAT(1), any
local geodesic of length $\pi$ is a geodesic. Thus, diam$(X) \geq
\pi$. Suppose there exists a geodesic $\gamma \co  [0,d] \to X$ with
$d>\pi$. Let $y=\gamma(0)$ and let $x=\gamma(\pi)$. Let $v \in \lk(x,X)$
be the outgoing tangent vector to $\gamma$ (ie, the tangent vector
to $\gamma \vert_{[\pi ,d]}$). By Lemma 3.3, there is a geodesic
$\alpha$ from $x$ to $y$ with tangent vector equal to $v$. But this
means that $\alpha$ and $\gamma \vert_{[\pi, d]}$ agree in a
neighborhood of $x$. This is clearly impossible since the distance
from $y$ decreases along $\alpha$ and increases along $\gamma
\vert_{[\pi, d]}$.\qed\enddemo

\proclaim{Lemma 3.5} If $A \in \a$, then $A$ is a Coxeter complex
$\sig$ or a suspension $\s^k * \sig$. 
\endproclaim

\demo{Proof}  Let $A \in \a$. We first show that the $(n-1)$--skeleton
of $A$, $A^{n-1}$, is a union of geodesic $(n-1)$--spheres which is
closed under reflection across each such sphere.

Suppose $x$ is a point in the relative interior of a
$(k-1)$--simplex $\sigma \subset A^{n-1}$. Then $\lk(\sigma,A)$ is
isometric to $\s^{n-k}$, hence it is an apartment in the
$(n-k)$--dimensional building $\lk(\sigma, X)$. The $(n-k-1)$--skeleton
of this apartment is a union of geodesic $(n-k-1)$--spheres closed
under reflection. Taking the join with $\s^{k-2} = \lk(x,\sigma)$, we see
that $\lk(x,A^{n-1})$ is a union of geodesic $(n-2)$--spheres in
$\lk(x,A)$ closed under reflections. It follows that, in a conelike
neighborhood of $x$ in $A$, the $(n-1)$--skeleton consists of a union
of geodesic $(n-1)$--disks, $D_1, D_2, \dots ,D_k$. In particular, any
geodesic through $x$ which enters this neighborhood through $A^{n-1}$
must also leave through $A^{n-1}$. Since this is true at every point
$x \in A$, we conclude that any geodesic in $A$ containing a
non-trivial segment in $A^{n-1}$, lies entirely in $A^{n-1}$.

Now let $D_i$ be one of the geodesic disks at $x$ as above. The
geodesic segments in $D_i$ emanating from $x$ extend to form a
geodesic $(n-1)$--sphere $H_i$ in $A$ which, by the discussion above, 
lies entirely in
$A^{n-1}$. We call $H_i$ a ``wall" through $x$. Reflection of $A$
across $H_i$ fixes $x$ and permutes the disks $D_1, \dots ,D_k$,
hence it permutes the walls through $x$. Moreover, if $H'$ is a wall
through some other point $x' \in A^{n-1}$, then $H'$ must intersect
$H_i$ (since they are two geodesic $(n-1)$--spheres in a $n$--sphere).
Say $z \in H_i \cap H'$. Applying the argument above with $x$
replaced by $z$ shows that reflection across $H_i$ takes $H'$ to some
other wall through $z$. Thus, it preserves $A^{n-1}$.

Let $W$ be the group generated by reflection across the walls of
$A$. By the previous lemma, $W$ acts on $A$ as a group of simplicial
isomorphisms. Since $A$ is a finite complex, $W$ is a finite
reflection group, or in other words, a spherical Coxeter group. Let
$A^W$ be the fixed set of $W$ which consists of the intersection of
all the walls.  Then $A^W$ is a geodesic $k$--sphere for some $k$, and
$A$ decomposes as a join, $A=A^W * \sig$.
\qed\enddemo

It follows from Lemma 3.5, that cells in $X$ are simplices or
suspensions of simplices.

\proclaim{Lemma 3.6} Any two cells $\sigma_1, \sigma_2$ in $X$
are contained in some $A \in \a$. \endproclaim

\demo{Proof} Since any cell is contained in an $n$--cell, it
suffices to prove the lemma for two $n$--cells. Let $x_1,x_2$ be
points in the interior of $\sigma_1,\sigma_2$, respectively. Let $\g$
be a geodesic from $x_1$ to $x_2$ and continue $\g$ to a geodesic of
length $\pi$. Let $y$ be the endpoint of $\g$ antipodal to $x_1$ and
note that
$\lk(x_1,X) = \lk(x_1, \sigma_1) \cong \s^{n-1}$. It follows from
Lemma 3.3 that the union of all geodesics from $x_1$ to $y$ forms a
subspace $A \in \a$. Since $A$ is a subcomplex and contains both
$x_1$ and $x_2$, it must contain $\sigma_1$ and $\sigma_2$.
\qed\enddemo

\proclaim{Lemma 3.7} If $A_1,A_2 \in \a$ share a common chamber
$\sigma$, then there is a simplicial isomorphism $\phi \co  A_1 \to A_2$
fixing $A_1 \cap A_2$ pointwise. Moreover, $A_1^W = A_2^W$.
\endproclaim

\demo{Proof} Let $x$ be a point in the relative interior of $\sigma$.
Let $p$ denote the north pole of $\s^n$ and fix an isometry 
$\theta \co  \lk(x,\sigma) \to \s^{n-1} = \lk(p,\s^n).$
Then there is a unique
isometry $\phi_i \co  A_i \to \s^n$ with $\phi_i (x)=p$ and the induced
map on $\lk(x,A_i)=\lk(x,\sigma)$ equal to $\theta$. Let $\phi =
\phi_2^{-1} \circ \phi_1$. Since the cell structure of $A_i$ is
completely determined by reflection in the walls of $\sigma$, the
isometry $\phi$ is also a simplicial isomorphism. For any point $y \in
A_1 \cap A_2$ not antipodal to $x$, there is a unique geodesic $\g$
from $x$ to $y$ which necessarily lies in $A_1 \cap A_2$. Since
$\phi_1$ and $\phi_2$ agree on the tangent vector $\g_x$, they agree
on all of $\g$. 

The last statement of the lemma follows from the fact that
$A_i^W=\sigma^W$. 
\qed\enddemo

It follows from Lemma 3.7, that $X$ itself decomposes as a join of
$A^W$ and a spherical building with Coxeter group $W$. This completes
the proof of Theorem 3.1.
\medskip

If we are not given an a priori cell structure, we can work in the
setting of locally spherical spaces and use the singular set (ie,
the branch set) of $X$ to define a cell structure. In this setting we
get the following theorem.

\proclaim{Theorem 3.8} Suppose $X$ is a locally spherical space of
dimension $n \geq 2$ satisfying
\roster
\item $X$ is CAT(1),
\item $\lk (x,X)$ is isometric to a building for every point $x \in X$.
\endroster
Then $X$ is isometric to a spherical building. The cell structure
determined by the singular set is that of a thick, spherical building
or a suspension of a thick, spherical building.
\endproclaim

\head More on spherical buildings \endhead

In contrast to the locally Euclidean case, a locally spherical
space which is simply connected and locally CAT(1) need not be
globally CAT(1). However, as we now show, under the stronger
hypothesis that links are isometric to buildings, a simply
connected locally spherical space of dimension $\geq3$ is CAT(1), 
and hence is also isometric to a spherical building.

\proclaim{Theorem 4.1} Suppose $X$ is a locally spherical space of
dimension $n \geq 3$ satisfying
\roster
\item $X$ is simply connected,
\item $\lk (x,X)$ is isometric to a building for every point $x \in X$.
\endroster
Then $X$ is CAT(1), hence it is isometric to a building.
\endproclaim

The hypothesis that $n \geq 3$ is essential here. In the 1980's
there was much interest in the the relation between 
incidence geometries and buildings (see for example \cite{T} and
\cite{K}). In \cite{T}, Tits proves a theorem analogous
to Theorem 4.1 for incidence geometries, with the same dimension
hypothesis. A counterexample in dimension $n=2$ is given by Neumaier in
\cite{N}. It is a finite incidence geometry of type $C_3$, with a 
transitive action of $A_7$ (the alternating group on 7 letters).
The flag complex
associated to Neumaier's $A_7$--incidence geometry is a 2--dimensional
simplicial complex, all of whose links are buildings, but which cannot
be covered by a building. One can metrize Neumaier's example by
assigning every 2--simplex the metric of a spherical
triangle with angles $\frac{\pi}{2}, \frac{\pi}{3}, \frac{\pi}{4}$. 
Passing to the universal
cover  gives a counterexample to Theorem 4.1 in dimension $n=2$.

Before proving the theorem, we will need some preliminary
lemmas. We begin with an easy consequence of Lemma 3.2.

\proclaim{Lemma 4.2} Let $X$ be as in Theorem 4.1 and
let $\g$ be a local geodesic in $X$ of length $\pi$ from $x$ to $y$.
Then there is a unique locally isometric extension $F_{\g} \co  \Sigma
(L_x) \to X$ of $\g$ (in the sense of Lemma 3.2, (1) and (2)).
\endproclaim

\demo{Proof} Suppose $A$ is an apartment in $L_x$ containing the tangent
vector $v=\g_x$. Then by Lemma 3.2, there is a
neighborhood $U$ of $v$ in $A$ and a unique locally isometric map $F_U
\co  \Sigma (U) \to X$ whose restriction to $\Sigma (v)$ is $\g$.
Using the maps $F_U$, we can extend $\g$ uniquely along any
geodesic in $L_x$ beginning at $v$. Since $L_x$ is simply connected
for $n \geq 3$, these extensions are compatible.
\qed\enddemo

Suppose $\g_1$ and $\g_2$ are two local geodesics of length $\pi$ from
$x$ to $y$ and let $v_i=(\g_i)_x$. Then it follows from the
construction of $F_{\g_i}$ that the following are equivalent.
\roster
\item $F_{\g_1}=F_{\g_2}$.
\item $\g_2$ is the restriction of $F_{\g_1}$ to $\Sigma(v_2)$ (and
vice versa).
\item There exists a locally isometric map of a spherical lune into $X$
with sides $\g_1$ and $\g_2$.
\endroster

We say that a geodesic $\eta \co  [0,a] \to X$ from $x$ is {\it
nonbranching} if any other geodesic $\eta'\co  [0,a] \to X$ from $x$ with
$\eta_x = \eta'_x$ is equal to $\eta$. (Or in other words, $\eta$ has
unique continuation at every point in its interior.) In particular,
if $\eta$ is contained in a cone-like neighborhood of $x$, then it is
non-branching.

The following is an immediate consequence of Lemma 4.2.

\proclaim{Lemma 4.3} Suppose $\g$ and $\eta$ are geodesics of length
$\leq \pi$ starting at $x$ and assume $\eta$ is nonbranching. Then
there is a locally isometric map of a spherical triangle into $X$
(possibly a geodesic or a spherical lune)
which restricts on two sides to $\g$ and $\eta$.
\endproclaim

The local isometry in the corollary above is essentially unique. More
precisely, we have the following.

\proclaim{Lemma 4.4} Let $T_1$ and $T_2$ be spherical triangles and 
$\Theta_1 \co  T_1 \to X$ and $\Theta_2 \co  T_2 \to X$ be
local isometries. If $\Theta_1$ and $\Theta_2$ agree along two edges
of the triangle, then one of the following holds.
\roster
\item $T_1=T_2$ (ie, they are isometric) and $\Theta_1=\Theta_2$,
\item $T_1$ and $T_2$ are hemispheres and $\Theta_1, \Theta_2$ agree
along their entire boundary.
\item $T_1$ and $T_2$ are spherical lunes and the two edges along
which they agree form one entire side of the lune.
\endroster\endproclaim

\demo{Proof} By hypothesis, $\Theta_1$ and $\Theta_2$ restrict along
two edges to local geodesics $\g$ and $\eta$ emanating from some point
$x$. The angle between these two edges is the distance in $L_x$
between $\g_x$ and $\eta_x$. Suppose this angle less than $\pi$.
Then clearly $T_1=T_2$. Since $X$ is locally CAT(1), the subspace of
$T_1$ on which $\Theta_1=\Theta_2$ must be locally convex, and hence
must be all of $T_1$.

If the angle between the edges is exactly $\pi$, then $T_i$ is either
a geodesic (if $\text{length}\,(\g) + \text{length}\,(\eta) < \pi$), a
spherical lune (if $\text{length}\,(\g) + \text{length}\,(\eta) =
\pi$), or a hemisphere (if $\text{length}\,(\g) +
\text{length}\,(\eta) > \pi$). In the last case, we may assume without
loss of generality that length$\,(\g)=\pi$. then $\Theta_1$ and
$\Theta_2$ are both restrictions of $F_{\g}\co \Sigma L_x \to X$. In
particular, they agree along $(\Sigma \g_x) \cup (\Sigma \eta_x)$ 
which forms the boundary of $T_i$.
\qed\enddemo

\proclaim{Lemma 4.5} The diameter of $X$ is at most $\pi$.
\endproclaim

\demo{Proof} The proof is the same as that of Lemma 3.4 (using
Lemma 4.2 in place of Lemma 3.3).
\qed\enddemo

\demo{Proof of Theorem 4.1} 
Fix a point $x$ in $X$. Define an equivalence relation on the set of
geodesics of length $\pi$ starting at $x$ by
$$\g_1 \sim \g_2 \Longleftrightarrow F_{\g_1}=F_{\g_2}$$
To prove Theorem 4.1, we define a covering
space $f \co  \tilde X \to X$ as follows. As a set, $\tilde X$ is defined as the
quotient 
$$\tilde X = \{\g \mid \text{$\g$ is a local geodesic of length $\leq \pi$
with $\g(0)= x$}\}/\sim$$
Note that only local geodesics of length $\pi$ can be identified in
$\tilde X$. 

The topology on $\tilde X$ is defined as follows.
Let $B_r(y)$ denote the ball of radius $r$ in $X$ centered at $y$.
Given a local geodesic $\g$ from $x$ to $y$ and a real number $r$ such
that $B_r(y)$ is conelike, define
$$\align
B_r(\g)=\{\eta \in \tilde X \mid &\text{ $\exists$ a locally isometric map of
a spherical triangle}\\
&\text{into $X$ which restricts on two sides to $\g$ and $\eta$}\\
&\text{and whose third side lies in $B_r(y)$}\}.
\endalign$$
If $\g$ has length $\pi$, then these locally isometric maps are all
restrictions of $F_{\g}$. In particular, $B_r(\g)$ depends only on the
class of $\g$ in $\tilde X$.
These sets form a basis for the topology on $\tilde X$. They also
define a metric (locally) on $\tilde X$. Namely, the distance between
$\g$ and $\eta \in B_r(\g)$ is the length of the third side of the
triangle. 

Define $f \co  \tilde X \to X$ to be the map taking $\g$ to its
endpoint. By Lemma 4.3, $f$ restricts to an isometry of $B_r(\g)$ onto
$B_r(y)$.
Letting $\g$ run over all local geodesics from $x$ to $y$, we claim
that these balls make up the entire inverse image of $B_r(y)$.
For suppose $\eta \in \tilde X$ with $f(\eta )=z \in
B_r(y)$. Since $z$ lies in a conelike neighborhood of $y$, the
geodesic $\delta$ from $z$ to $y$ is nonbranching. It follows from
Lemma 4.3, applied to $\eta$ and $\delta$, that there exists a
local isometry of a spherical triangle into $X$ which restricts on two
side to $\eta$ and $\delta$. The restriction to the third side, $\g$,
is a local geodesic from $x$ to $y$ such that $\eta \in B_r(\g)$. 

It remains to show that for distinct $\g_i \in \tilde X$, the balls
$B_r(\g_i)$ are disjoint. Suppose $\delta \in B_r(\g_1) \cap
B_r(\g_2)$. Then there is a local isometry $\Theta_i$ of a spherical
triangle with sides $\g_i$ and $\delta$, for $i=1,2$. Since $B_r(y)$ is
conelike, there is a unique geodesic $\eta$ from $y$ to the endpoint
$z$ of $\delta$. Thus $\Theta_1$ and
$\Theta_2$ agree along two edges, $\delta$ and $\eta$. It follows from 
Lemma 4.4 that $\g_1=\g_2$ in $\tilde X$.

This proves that $f\co \tilde X \to X$ is a covering map. By hypothesis,
$X$ is simply connected, and it is easy to verify that $\tilde X$ 
is connected, thus $f$ is injective. In particular, for any $y
\in X$ of distance less than $\pi$ from $x$, there is a unique local 
geodesic from $x$ to $y$. Moreover, it follows from Lemma 4.3 that
this geodesic varies continuously with the endpoint $y$. Since $x$ was
chosen arbitrarily, this applies to all $x$ and $y$. By Theorem 2.2(2), we
conclude that $X$ is CAT(1).
\qed\enddemo

\head Euclidean buildings \endhead

\proclaim{Theorem 5.1} Suppose $X$ is a connected, locally
Euclidean complex satisfying 
\roster
\item $X$ is CAT(0),
\item for every point $x \in X$, $L_x = \lk(x,X)$ is isometric to a
spherical building.
\endroster
Then $X$ decomposes as an orthogonal product $X=\r^l \times X_1 \times
\dots \times X_k$, where $l \geq 0$, and each $X_i$ is one of the following,
\roster
\item a thick, irreducible Euclidean building,
\item the Euclidean cone on a thick, irreducible spherical building, 
\item a tree.
\endroster
\endproclaim

\remark{Remark} The reader may object that a tree is a 1--dimensional
irreducible Euclidean building whose apartments are Coxeter complexes
for the infinite dihedral group. However, the standard building metric
on a 1--dimensional Euclidean building would assign the same length to
every edge of the tree. Since this need not be the case in our situation,
we list these factors separately.
\endremark
\medskip

In \cite{KL}, Kleiner and Leeb introduce a more general notion of a
Euclidean building, which we will call a ``metric Euclidean
building'', and prove an analogous product decomposition theorem
(Prop. 4.9.2) for these buildings. We review their definition in the
context of locally Euclidean spaces. (Kleiner and Leeb work in a
more general setting.) 

Call a
group $W$ of affine transformations of $\Bbb R^n$ an {\it affine Weyl
group} if $W$ is generated by reflections and the induced group of
isometries on the sphere at infinity is finite. Affine Weyl groups
include Euclidean and spherical Coxeter groups, as well as nondiscrete
groups generated by reflections across parallel walls.

Let $\a$ be a collection of isometric embeddings of $\Bbb R^n$ into a
locally Euclidean space $X$ of dimension $n$. We call $\a$ an {\it
atlas} for $X$ and the images of the embeddings are called {\it apartments}.

\definition{Definition 5.2} Suppose $X$ is a CAT(0), locally Euclidean
space of dimension $n$. Then $X$ is a {\it metric Euclidean Building}
if there is an atlas $\a$ and an affine Weyl group $W$ such that
\roster
\item Every geodesic segment, ray, and line is contained in an apartment.
\item $\a$ is closed under precomposition with $W$.
\item If two apartments $\phi_2(\Bbb R^n), \phi_2(\Bbb R^n)$
intersect, then $\phi_1^{-1} \circ \phi_2$ is the restriction of some
element of $W$.
\endroster
\enddefinition 

(In the context of locally Euclidean spaces, this definition agrees
with that of Kleiner and Leeb since their first two axioms hold
automatically for locally Euclidean spaces.) It is an immediate
consequence of Theorem 5.1 that the space $X$ is a metric Euclidean
building.

\proclaim{Corollary 5.3} Let $X$ be as in Theorem 5.1. Then $X$ is a
metric Euclidean building.
\endproclaim

 Conversely, it is easy to see that a metric Euclidean
building satisfies the hypotheses of Theorem 5.1. Thus, for locally
Euclidean spaces, Theorem 5.1 also provides another proof of
Kleiner and Leeb's product decomposition theorem.
\medskip

The proof of Theorem 5.1 will occupy the remainder of this
section. As in the spherical case, the key is to find lots of apartments. 
By Proposition 2.3, the apartments in $X$ are isometrically embedded
copies of $\Bbb R^n$, known as {\it $n$--flats}. 
The crucial step to constructing $n$--flats is to find {\it flat
strips} (isometrically embedded copies of $\r \times [0,a]$ for some
$a>0$).

\definition{Definition 5.4} 
Let $X$ be a $CAT(0)$--space.  We will call the triangle $\Delta$ in
$X$ a {\it Euclidean (or flat) triangle}, if its convex hull is
isometric to a triangle in $\Bbb R^n$. 
\enddefinition

For a triangle $\Delta (x,y,z)$ we denote the segment from $x$ to $y$
by $xy$, etc. The angle between $xy$ and $xz$ is defined as the
distance in $L_x$ between the tangent vectors to $xy$ and $xz$
and it is denoted by $\angle_x (xy, xz)$.
The following lemma follows immediately from Proposition 3.13 of
\cite{Ba}.

\proclaim{Lemma 5.5}
 Let $X$ be a $CAT(0)$ space, $\Delta = \Delta(a,b,c)$ a triangle 
in $X$. Let $d$ be a point between  $a$ and $b$. Suppose the triangles 
$\Delta(a,d,c)$ and $\Delta(b,d,c)$ are Euclidean.
 If in addition  $\angle _d (da,dc) + \angle _d (db,dc) = \pi$,
then the original triangle  $\Delta(a,b,c)$ is Euclidean.
\endproclaim
   
 The condition on the angles is automatically satisfied if the
 geodesic $ab$
is non-branching, for example if $b$ lies in a cone-like
neighborhood of $a$. 
 
If $L$ is a locally spherical space and $r>0$, let $C_r(L)$ denote the
Euclidean cone on $L$ of radius $r$ (ie, the geodesics emanating from
the cone point have length $r$).
The following lemma is an analogue of Lemma 3.2. The proof is
essentially the same and the details are left to the reader.

\proclaim{Lemma 5.6} Suppose $X$ is a locally Euclidean space of
dimension $n \geq 2$ such that the link of every point in $X$ is
isometric to a spherical building. Let $\g$ be a locally geodesic ray
from $x$ and let $A_x$ be an apartment in $L_x$ containing
$\g_x$. Then for any $r>0$, there is a neighborhood $N_x \subset A_x$ 
of $\g_x$ and a unique locally isometric map $\Theta$ of the Euclidean
cone $C_r(N_x)$ into $X$ such that
\roster
\item for any $v \in N_x$, the restriction of $\Theta$ to $C_r (v)$ is a
local geodesic with tangent vector $v$, and
\item the restriction of $\Theta$ to $C_r (\g_x)$ is precisely $\g 
\vert_{[0,r]}$.
\endroster
\endproclaim

Suppose $\g$ is as in Lemma 5.6, and $\alpha$ is a geodesic in $L_x$
originating at $\g_x$. Then Lemma 5.6 implies that there is a unique
extension of $\g \vert_{[0,r]}$
to a locally isometric map of $C_r(\alpha)$ into $X$. If $X$ is
CAT(0), the map is an isometric embedding. This enables us to construct
Euclidean triangles in $X$ since the image of any triangle in
$C_r(\alpha)$ is Euclidean.

From now on, we assume that $X$ satisfies the hypotheses of Theorem
5.1. 

\proclaim{Lemma 5.7} Let $\gamma\co \r \to X$ be a geodesic ray 
with $x=\g(0)$. If $y$ lies in a conelike neighborhood of $x$, then
for any $t \in \r$, the triangle $\Delta(x,y,\gamma(t))$ is Euclidean.
\endproclaim 

\demo{Proof} Assume, without loss of generality, that $t>0.$ Since $y$ 
lies in a conelike neighborhood of $x$, the geodesic $\eta$ from $x$
to $y$ is non-branching. Choose a geodesic $\alpha$ in $L_x$ from
$\g_x$ to $\eta_x$ and extend $\g \vert_{[0,t]}$ to an isometric
embedding $\Theta$ of $C_t(\alpha)$ into $X$. Since $\eta$ is
non-branching, it agrees with $\Theta$ on $C_t(\eta_x)$. Thus,
$x,y,\gamma(t)$ span a Euclidean triangle.
\qed\enddemo

\proclaim{Lemma 5.8} Let $\g$, $x$, and $y$ be as above and let $\eta$
 be the geodesic from $x$ to $y$.
 Suppose $\angle_x (\eta ,\gamma ^+) + \angle_x (\eta ,\gamma ^-) =
 \pi$. Then $y$ and $\gamma$ span a Euclidean strip. 
\endproclaim

\demo{Proof} By the previous lemma, for each $t \in \Bbb R$, the
 triangle $\triangle (x,y,\gamma(t))$ is Euclidean. By Lemma 5.5,
 every triangle of the form $\triangle (y, \g(t_1), \g(t_2))$ is
 Euclidean. Since any two points in the span of $y$ and $\g$ lie on
 such a triangle, the lemma follows.
\qed\enddemo

For any subspace $Y$ of $X$, and any point $x \in Y$, we
denote the link of $x$ in $Y$ by $L_x Y$.

\proclaim{Lemma 5.9} Let $F$ be an $m$--flat in $X$. Let $\eta \co  [0,r] 
\to X$ be a geodesic from $x \in F$ to a point $y$ lying in a
 conelike neighborhood of $x$. Suppose that the distance in $L_x$
 from $\eta_x$ to any point in $L_x F$ is $\frac{\pi}{2}$. Then $y$ and $F$ 
 span a flat $R^m \times[0,r]$.
\endproclaim

\demo{Proof} Let $Z=\Bbb R^m \times [0,r]=F \times [0,r]$. Let $Y$ be 
the subspace of $X$ spanned by $F$ and $y$.  Consider the natural map
$f \co  Z \to Y$ that takes $F \times \{0\}$ via the identity map to $F$
and takes the line segment between $(z,0)$ and $(0,r)$ to the geodesic
in $Y$ from $z$ to $y$. By Lemma 5.8, the restriction of $f$ to the
strip $\g \times [0,r]$ is an isometry for any geodesic $\g \co  \Bbb R
\to F$ through $x$.  We must show that $f$ is isometric on all of $Z$.

Any two points $y_1, y_2$ in $Y$ lie on the image of a
triangle $T \subset Z$ with vertices $(0,r), (z_1,0), (z_2,0)$. By the
discussion above, $T$ is a comparison triangle for its image $f(T)$ in
$Y$. Hence, by the CAT(0) condition, the distance between $y_1$ and
$y_2$ is at most the distance between the corresponding points in
$T$. Thus, $f$ is distance non-increasing.
Moreover, if $y_1$ or $y_2$ lies on $\eta$, then they both lie in
a Euclidean strip, as above, so these two distances agree.

To prove the reverse inequality, choose $r' < r$ and let $y'
=\eta(r')=f((0,r'))$ (Figure 1). Consider the induced map $df$
between the links $L_{(0,r')} Z \cong \s^m$ and $L_{y'} Y$.  It
suffices to prove that $df$ is an isometry, for in this case, the
triangle with vertices $y', y_1, y_2$ has the same angle and same two
side lengths at $y'$ as its comparison triangle in $Z$,
so by the CAT(0) condition, the opposite side is at least as
long as in the comparison triangle.

To see that $df$ is an isometry, note that the fact that $f$ restricts
to an isometry on strips $\g \times [0,r]$ implies that $df$ maps
points of distance $\pi$ in $L_{(0,r')} Z$ to points at distance
$\pi$ in $L_{y'} Y$. On the other hand, since $f$ is distance
non-increasing at all points, it must also be distance non-increasing
on links. But these two facts contradict each other unless $df$ is an
isometry.
\qed\enddemo

\figure
\relabelbox\small
\epsfysize=3truein
\epsffile{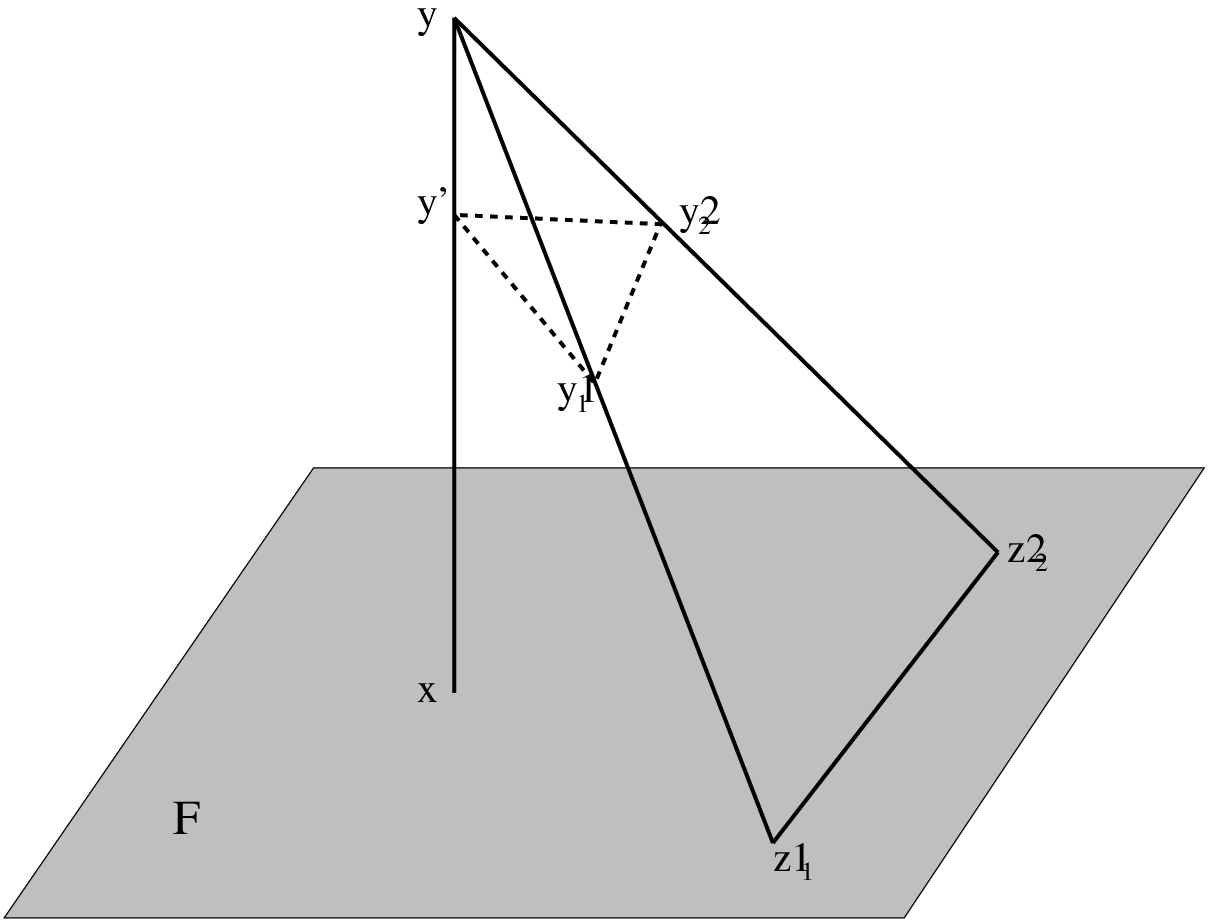}
\relabel {y}{$y$}
\adjustrelabel <-1pt,-1pt> {y1}{$y_1$}
\relabel {y2}{$y_2$}
\adjustrelabel <-1pt,-1pt> {y'}{$y'$}
\relabel {z}{$z$}
\relabel {z1}{$z_1$}
\relabel {z2}{$z_2$}
\relabel {x}{$x$}
\relabel {F}{$F$}
\endrelabelbox
\endfigure

\proclaim{Lemma 5.10} Every geodesic and every flat strip in $X$ is
contained in an $n$--flat.
\endproclaim

\demo{Proof}
 Let $F$ be an $m$--flat in $X$ with $m<n$. 
Let $a \in [-\infty ,0]$, $b \in [0,\infty ]$
be chosen so that $F$ is contained in a flat embedded  $F\times[a,b]$
(with $F = F \times \{0\}$) and such that $a,b$ are maximal, (that is, 
one cannot embed this  $F\times[a,b]$ in a bigger
$F\times[a_0,b_0]$).
 
 We claim that $a=-\infty ,b =\infty $ . Assume the contrary.  Say,
for example, $b <\infty $. Let $x$ be a point in $F\times \{b\}$.
$F\times[a,b]$ determines an $m$--dimensional hemisphere $H$ in $L_x X$
containing $L_x F$.  (If $a=b=0$, choose any hemisphere containing
$L_x F$.)  Choose an $m$--sphere in $L_x X$ containing $H$ and let $v $ be
the point of the sphere, which has distance $\frac{\pi}{2}$ from $H$. Let
$\eta$ be a ray emanating from $x$ in the direction of $v$. By the
previous lemma, for small $r$, $\eta (r)$ and $F\times
\{b\}$ span a flat $F \times [0,r]$.  The choice of $v$, together
with Lemma 5.5, guarantee that $F\times[a,b]$ and $F\times[0,r]$
fit together to form a flat strip $F\times[a,b+r]$. This
contradicts the maximality of $b$.

Thus, any $m$--flat, and more generally, any geodesically embedded
$\r^m \times [a,b]$, $m<n$, can be embedded in an $(m+1)$--flat. It
follows inductively, that every geodesic and every flat strip is
contained in an $n$--flat.
\qed\enddemo

\demo {Proof of Theorem 5.1} Let $F$ be an $n$--flat in $X$. The set
$Y$ of the singular points in $F$ is closed and locally it is a union
of hyperplanes, so $Y$ is globally a union of hyperplanes too. We call
these the singular hyperplanes.  The set of singular hyperplanes is
locally invariant under reflection in each of these hyperplanes, so if
two singular hyperplanes $H_1$ and $H_2$ are not parallel, the
reflection, $H'_1$, of $H_1$ across $H_2$ is also singular. Moreover,
if $H_1$ and $H_2$ are parallel, and there exists a singular
hyperplane not orthogonal to $H_1$ and $H_2$, then a simple exercise
shows that $H'_1$ can be obtained by a series of reflections across
intersecting singular hyperplanes. Thus, again, $H'_1$ must be a
singular hyperplane. 

Let $Y_1,Y_2,\ldots,Y_k$ be a maximal decomposition of $Y$ into
mutually orthogonal families of singular hyperplanes. It follows from
the discussion above that each $Y_i$
is either a discrete family of parallel hyperplanes or the set of
walls of an irreducible spherical or Euclidean Coxeter group (compare
\cite{Br}, VI.1).

Taking $F_i$ to be the subspace of $F$ generated by
the normal vectors to the hyperplanes in $Y_i$, we obtain an
orthogonal decomposition of $F$ into $F=F_0 \times F_1 \times \cdots \times
F_k$, where $F_0$ is a subspace parallel to all
of the singular hyperplanes. If $Y_i$ is a family of parallel
hyperplanes, then $F_i$ is $1$--dimensional. Otherwise, $Y_i$ gives
$F_i$ the structure of a Euclidean Coxeter complex or an (infinite)
Euclidean cone on a spherical Coxeter complex.
 
We now choose an $n$--flat $E$ such that, in the decomposition $E=E_0
\times E_1 \times \cdots \times E_k$, the dimension $l$ of the
factor $E_0$ is minimal, and among all flats with dim$(E_0)=l$, we
require that $E$ minimize the number $k$. (As we will see below, the
number $k$ is actually independent of the choice of $E$.)  In this
flat $E$, we choose a point $x$, such that for $j=1,\ldots,k$, the
hyperplanes $Y_j$ determine a full-dimensional, irreducible, spherical
Coxeter group in the sphere $L_x E_j$. Intuitively, we have chosen $E$
to be the most complicated $n$--flat and $x$ to be the most complicated
point in $E$.

The product decomposition of $E$ gives rise to a decomposition of $L_x
E$ as an orthogonal join, $L_x E =  L_x E_0 * L_x E_1 * \cdots * L_x
E_k $. Since $L_x E$ is an apartment in the spherical building $L_x
X$, it follows that $L_x X$ splits as a 
join of $k$ irreducible buildings and an $(l-1)$--sphere,
$L_xX= L_0*L_1*\cdots *L_k$.

 Let $G$ be another flat through $x$. As above, we can factor $G$ as
an orthogonal product $G = G_0 \times G_1 \times \cdots \times G_{k'}$
Since $L_x G$ is also an apartment in $L_x X$, it must have the same
simplicial structure as $L_x E$. Clearly this is possible only if
$\text{dim}\,(G_0) = \text{dim}\,(E_0)$, $k'=k$ and (up to
permutation) $L_x G_j \cong L_x E_j$ for all $j$. (Note that although
we have shown that the simplicial structure of $G_j$ and $E_j$ agree
in a neighborhood of $x$, we do not yet know that they agree
globally.)

If $y$ is any other point in $X$, then by Lemma 5.10, there is a flat
$G$ containing $x$ and $y$, so the decomposition of $G$ gives rise to
a decomposition of $L_yX$ as a join $L_{y,0} * L_{y,1}* \cdots *
L_{y,k}$ of spherical buildings.  In this case, however, the factors
need not be irreducible. (Consider, for example, the case where $y$ is
a nonsingular point.)

Next, we prove that if the decomposition of $E$ is not the trivial
one, we can decompose $X$ as a product.  Let $X_i$ be the union of the
geodesic rays emanating from $x$ in the direction of $L_i$. Then $X_i$
is connected and by Lemma 5.7, it is locally convex, hence it is
convex.  To prove that $X = X_0\times X_1 \times \cdots \times X_k$,  
we apply Theorem II.9.24 of \cite{BH} which states that splittings
of $X$ as a product correspond to splittings of the Tits boundary
$\partial_TX$ as a join. The Tits boundary may be viewed as the set of
rays emanating from $x$. To show that $X = X_0\times X_1 \times
\cdots \times X_k$, it suffices to show that $\partial_TX =
\partial_TX_0*\partial_TX_1 * \cdots * \partial_TX_k$. For this we
must verify (cf~\cite{BH} Lemma II.9.25):
\roster
\item Every ray $\g$ lies in the span of rays $\g_0,\g_1, \dots ,\g_k$
with $\g_i \subset X_i$.
\item Rays $\g_i, \g_j$ lying in distinct factors $X_i, X_j$ have
distance $\frac{\pi}{2}$ in the Tits metric.
\endroster

The first of these conditions follows from the fact that any ray $\g$
at $x$ is contained in an $n$--flat $G$ which, by the discussion above,
decomposes as a product $G=G_0 \times G_1 \times \dots \times G_k$
with $G_i \subset X_i$. 

For the second condition, it is enough to
show, that any geodesic $\gamma$ in $X_i$ through $x$ and any ray
$\eta$ in $X_j$ emanating from $x$ ($i\neq j$) span a flat halfplane.
If not, then there is a largest $t$ such that $\eta(t)$ and
$\gamma$ span a flat strip. By Lemma 5.8, we can assume that $t$ is
bigger than $0$.  By Lemma 5.10, we can embed this strip in an
$n$--flat $T$.  Let $z=\eta(t)$ and let $\gamma _1$ be the line in $T$
parallel to $\gamma$ through $z$.  The product decomposition of $T$
gives rise to a join decomposition of the building $L_z X = L_{z,0}* 
 L_{z,1} * \cdots * L_{z,k}$.  The tangent vectors to $\gamma_1$ at $z$
lie in $L_{z,i}$ while the tangent vectors to $\eta$ at $z$ lie in
$L_{z,j}$. It follows that the angle condition of Lemma 5.8 is
satisfied, so that for small $\epsilon$, $\eta(\epsilon +t)$ and
$\gamma _1$ span a flat strip. This strip fits together with the strip
between $\gamma _1$ and $\gamma$, contradicting the maximality of $t$.
This proves condition (2) and we have shown that $X$ decomposes as an
orthogonal product of the $X_i$'s.

It remains to identify the factors $X_i$. Since $X_i$ is convex
in $X$, it is CAT(0), and at any point $y \in X_i$, $L_yX_i=L_{y,i}$ 
is a building. Thus, each factor satisfies the hypotheses of Theorem 5.1.
 The factor $X_0$ contains no singular points, for if
$y \in X_0$ is singular, then we could find an $n$--flat $F$ through
$y$ with $dim F_0 < dim E_0$. In other words, $X_0=\r^l$. 

Assume $i>0$.  If $X_i$ is 1--dimensional, then the CAT(0) condition
implies that it is a tree. If $\text{dim}\,X_i \geq 2$, then $X_i$
contains an apartment $E_i$ which is either an irreducible Euclidean
Coxeter complex or the cone on an irreducible spherical Coxeter
complex. Thus, chambers in $E_i$ are simplices or cones on simplices,
and the simplicial structure of $E_i$ is completely determined by a
single chamber. It follows that any other apartment in $X_i$ sharing a
chamber with $E_i$ has the same simplicial structure. By Lemma 5.10,
any two points (hence any two chambers) in $X_i$ lie in a common
apartment, so any two apartments are isomorphic. If $E_i$ is a
Euclidean Coxeter complex, then we conclude that $X_i$ is an
irreducible Euclidean building. In the case that $E_i$ is the cone on
a spherical Coxeter complex, it has only one vertex, namely $x$, which
must be the cone point for every apartment. Thus, in this case, $X_i$
is the cone on the irreducible spherical building $L_i$.  In either
case, since the simplicial structure on $X_i$ is defined by its
singular set, these buildings must be thick. This completes the proof
of Theorem 5.1.
\qed\enddemo

It is reasonable to ask whether an analogous theorem holds for
locally hyperbolic spaces; that is, if $X$ is locally
hyperbolic and CAT(-1), and every link in $X$ is isometric to a
spherical building, can we conclude that $X$ is a hyperbolic building
(ie, a buildings whose apartments are copies of $\Bbb H^n$
cellulated by the walls of a discrete, hyperbolic reflection group)?
The answer is no. Although one can construct lots of embedded copies
of $\Bbb H^n$ under these hypotheses, these ``apartments'' need not
have the structure of a Coxeter complex. For example, let $\g_1, \dots
\g_5$ be five geodesics $\g_i\co \r \to \Bbb H^2$ which form a right-angled
pentagon. The reflections across these lines generate a hyperbloic
Coxeter group. Now let $\g_6$ be another geodesic line intersecting
$\g_1$ orthogonally and disjoint from the other $\g_i$'s.  We can
choose $\g_6$ so that the three geodesics intersecting $\g_1$ generate
a nondiscrete group of reflections. Then $\g_1, \dots ,\g_6$ cannot be
the walls of a reflection group acting on $\Bbb H^2$.  (Note that
these six lines are locally, but not globally closed under
reflections.)  It is now possible to construct a simply connected,
locally hyperbolic space $X$, all of whose links are spherical
buildings (in fact we can take the link of every vertex to be the
$K_{3,3}$ graph with edge lengths $\frac{\pi}{2}$) and such that $X$ contains
an isometrically embedded $\Bbb H^2$ whose singular set consists of
the six lines above. This cannot be a hyperbolic building.

\head One-dimensional spherical buildings \endhead

In this section we give a metric characterization of one-dimensional
spherical buildings. An equivalent characterization appears as Lemma
6.1 in \cite{BB2}. We include a proof here for completeness.

\proclaim{Theorem 6.1} Suppose $X$ is a connected, one-dimensional
piecewise spherical complex satisfying
\roster
\item $X$ is CAT(1) and $\text{diam\,}(X)=\pi$,
\item every vertex of $X$ has valence $\geq 3$.
\endroster
Then $X$ is either a thick spherical building or $X=\Sigma Y$ for
some discrete set $Y$. \endproclaim

\remark{Remark} We can also apply this theorem in cases where every
vertex has valence $\geq 2$ simply by ignoring those vertices of
valence $2$. In this case, however, some information about the
original cell structure will be lost since vertices of valence $2$ are
invisible to the metric.
\endremark

Assume throughout this section that $X$ satisfies the hypotheses of
Theorem 6.1. In this situation, the analogue of Lemma 3.3 is easy to
prove. 

\proclaim{Lemma 6.2} Let $x,y \in X$ be antipodal points and $v \in
\lk(x,X)$ a tangent vector at $x$. Then there exits a unique geodesic
$\g$ from $x$ to $y$ such that $\g_x =v$. \endproclaim

\demo{Proof} Choose a point $z$ in the open star of $x$ so that the
geodesic $\alpha$ from $x$ to $z$ satisfies $\alpha_x = v$. We first
observe that $z$ cannot be antipodal to $y$. For if $w$ is a point on
$\alpha$ between $x$ and $z$, then the geodesic from $w$ to $y$ must
pass through either $z$ or $x$. If both $z$ and $x$ are antipodal to
$y$, then $d(w,y)>\pi$, contradicting $\text{diam\,}(X)=\pi$. Thus,
$d(x,z)<\pi$ and hence there is a unique geodesic $\beta$ from $z$ to
$y$. This geodesic cannot pass through $x$ (since this would imply
$d(z,y)>d(z,x)=\pi$), so it must leave $z$ in the opposite direction
from $\alpha$. Thus, the concatenation $\alpha \cdot \beta$ is a
local geodesic from $x$ to $y$. Moreover, since $\alpha$ contains no
antipodal points to $y$ (other than $x$), the length of $\alpha
\cdot \beta$ cannot be more than $\pi$. Thus, $\g = \alpha \cdot
\beta$ is the desired geodesic.
\qed\enddemo

\proclaim{Lemma 6.3} Suppose $x,y \in X$ are antipodal points. If $x$
is a vertex of $X$, then so is $y$. \endproclaim

\demo{Proof} If $x$ is a vertex of $X$, then $\lk (x,X)$ contains at
least 3 distinct points $v_1, v_2, v_3$. By the previous lemma, there
are geodesics $\g_1, \g_2, \g_3$ from $x$ to $y$ with
$(\g_i)_x=v_i$. Any two of these geodesics intersect only at $x$ and
$y$ since, otherwise, they would form a circuit of length $<2\pi$. 
Thus, $\lk (y,X)$ has at least 3 distinct points, so $y$ is a vertex.
\qed\enddemo

\proclaim{Lemma 6.4} All edges of $X$ have the same length, 
$\frac{\pi}{m}$ for some integer $m \geq 1$. \endproclaim

\demo{Proof} Since $X$ is connected, if any two edges have different
lengths, then we can find a pair of adjacent edges, $e_1,e_2$ of
different lengths. Let $x$ be their common vertex and let $x_1,x_2$
be the other endpoints of $e_1,e_2$ respectively. Say
$d(x,x_1)<d(x,x_2)$. Choose
a point $y$ antipodal to $x$ and a third edge $e_3$ emanating from
$x$. By Lemma 6.2, there are geodesics $\g_1,\g_2,\g_3$ from $x$ to
$y$ which begin with $e_1,e_2,e_3$ respectively. Any two of these
geodesics, $\g_i,\g_j$, form a loop of length $2\pi$. Consider the
loop $\g_1 \cdot \g_3$. Let $z$ be the point on $\g_3$ antipodal
$x_1$. By Lemma 6.3, $z$ is a vertex. Next, consider the loop $\g_2
\cdot \g_3$. Let $z'$ be the point on $\g_2$ antipodal to $z$. By
Lemma 6.3, $z'$ is a vertex. But $d(x,z')=d(y,z)=d(x,x_1)<d(x,x_2)$.
Thus, $z'$ lies in the interior of the edge $e_2$. This is a
contradiction. We conclude that all edges of $X$ have the same
length. Since antipodal points to vertices are also vertices, this
length must be $\frac{\pi}{m}$ for some integer $m \geq 1$.
\qed\enddemo

If all edges of $X$ have length $\pi$, then the distance between
any two vertices of $X$ is exactly $\pi$. It follows that $X$ has
only two vertices (otherwise, it would have diameter $>\pi$).
Thus, $X$ is the suspension of the discrete set $Y$
consisting of the midpoints of the edges.

From now on, we will assume that every edge has length $\frac{\pi}{m}$ for
some $m \geq 2$. In this case, every circuit of length $2\pi$ is a
Coxeter complex for the dihedral group $D_{2m}$. We will refer to
such circuits as apartments.

\proclaim{Lemma 6.5} Any two simplices of $X$ are contained in a common
apartment. If $A_1$ and $A_2$ are apartments containing a
common edge $e$, then there is a simplicial isomorphism $\phi \co  A_1
\to A_2$ fixing $A_1 \cap A_2$ pointwise. \endproclaim

\demo{Proof} It suffices to prove the first statement for two edges
$e_1,e_2$. Choose a point $x$ in the interior of $e_1$. Suppose $e_2$
contains a point $y$ antipodal to $x$ (which, by Lemma 6.3 must lie
in the interior of $e_2$). Then $\lk (x,e_1)$ consists of two
points $v_1,v_2$ each of which gives rise to a geodesic $\g_1,\g_2$
from $x$ to $y$. Together, these two geodesics form an apartment
containing $e_1$ and $e_2$.

Suppose, on the other hand, that no point of $e_2$ is antipodal to
$x$. Let $\g_1$ be a geodesic from $x$ to a point $y$ in the interior
of $e_2$. Extend $\g_1$ to a geodesic of length $\pi$. Note that this
extended geodesic (which we still denote $\g_1$) must contain all of
$e_2$. Let $z$ be the endpoint of $\g_1$, so $z$ is antipodal to 
$x$. Then there exists a geodesic $\g_2$ from $x$ to $z$ which begins
along $e_1$ in the direction opposite to $\g_1$. The union, $\g_1 \cup
\g_2$ is the desired apartment.

The second statement is obvious since the CAT(1) hypothesis implies
that $A_1 \cap A_2$ is connected.
\qed\enddemo

\head Main theorems \endhead

Combining the results of the previous sections, we arrive at our main
theorems. 

\proclaim{Theorem 7.1} Suppose $X$ is a connected, piecewise spherical
complex of dimension $n \geq2$ satisfying the following conditions.
\roster
\item $X$ is CAT(1).
\item Every $(n-1)$--cell is contained in at least two $n$-cells.
\item The link of every $k$--cell, $k \leq n-2$, is connected.
\item The link of every $(n-2)$--cell has diameter $\pi$.
\endroster
Then $X$ is isometric to a spherical building. The cell structure
determined by the singular set is that of a thick spherical building
or a suspension of a thick spherical building.
\endproclaim

\demo{Proof} Let $L_{\sigma}= \lk (\sigma,X)$. By Theorem 2.2, if $X$
is CAT(1), then so is $L_\sigma$ for every cell $\sigma$. In
particular, for an $(n-2)$--cell $\sigma$, $L_\sigma$ is a
1--dimensional, piecewise spherical complex which is CAT(1), diameter
$\pi$, and has every vertex of valence $\geq 2$. Ignoring vertices of
valence $2$ gives a complex satisfying Theorem 6.1. Thus, $L_\sigma$
is isometric to a 1-dimensional spherical building. The theorem now
follows by induction from Theorem 3.8.
\qed\enddemo

In the theorem above, we could have assumed that $X$ was locally
spherical instead of piecewise spherical if we interpret ``$k$--cell''
as meaning $k$--dimensional strata of the singular set.  In fact, in
ignoring vertices of valence $2$, some information about the given
cell structure on $X$ may be lost. For example, {\it any} cell
decomposition of the standard 2--sphere satisfies the conditions of the
theorem, but need not be the cell structure of a building. To
guarantee that the original cell structure is reflected in the metric,
we would need to assume that $X$ is thick.

\proclaim{Theorem 7.2} Suppose $X$ is a connected, piecewise spherical
complex of dimension $n \geq2$ satisfying the following conditions.
\roster
\item $X$ is CAT(1).
\item Every $(n-1)$--cell is contained in at least three $n$--cells.
\item The link of every $k$--cell, $k \leq n-2$, is connected.
\item The link of every $(n-2)$--cell has diameter $\pi$.
\endroster
Then, with respect to its given cell structure, $X$ is a spherical building.
\endproclaim

\demo{Proof} The proof is the same using Theorem 3.1 instead of
Theorem 3.8.
\qed\enddemo

Analogous results hold in the piecewise Euclidean setting.

\proclaim{Theorem 7.3} Suppose $X$ is a connected, piecewise Euclidean
complex of dimension $n \geq 2$ satisfying the following conditions.
\roster
\item $X$ is CAT(0).
\item Every $(n-1)$--cell is contained in at least two $n$--cells.
\item The link of every $k$--cell, $k \leq n-2$, is connected.
\item The link of every $(n-2)$--cell has diameter $\pi$.
\endroster
Then $X$ is a metric Euclidean building. If, in addition, $X$ is thick,
then $X$ is a product of irreducible Euclidean buildings and trees
with respect to its given cell structure.
\endproclaim

\demo{Proof} The first hypothesis implies that $\lk (\sigma,X)$ is
CAT(1) for every $\sigma$. By Theorem 7.1, we conclude that for $v$ a
vertex of $X$, $\lk (v,X)$ is isometric to a spherical building.
It follows by Theorem 5.1 and Corollary 5.3 that $X$ is a metric
Euclidean building and that it factors as a product
of irreducible Euclidean buildings, cones on spherical buildings,
trees, and a nonsingular Euclidean space. If $X$ is thick, the
components of the nonsingular set of $X$ are the interiors of the
$n$--cells. In particular, they are bounded. Thus, only the first and
last type of factor can occur in this situation.
\qed\enddemo

Replacing the CAT(1) condition in Theorem 7.1 by a simply connectedness
condition does not give a satisfying characterization because of
the problem in dimension 2. It does, however, give an analogue of
Tit's theorem about incidence geometries (\cite{T}, Theorem 1).								

\proclaim{Theorem 7.4} Suppose $X$ is a piecewise spherical
(respectively Euclidean)  complex of dimension $n\geq 3$ satisfying
\roster
\item $X$ is simply connected.
\item The link of every $k$--cell, $k \leq n-4$, is simply connected.
\item The link of every $(n-3)$--cell is isometric to a building.
\endroster
Then $X$ is isometric to a spherical building (respectively metric Euclidean
building).
\endproclaim

\demo{Proof} In the spherical case, the theorem follows from Theorem
4.1 and induction. In the Euclidean case, links in $X$ are spherical
buildings (by the spherical case of the theorem) so $X$ is locally
CAT(0). Since $X$ is simply connected, it is also globally
CAT(0) (Theorem 2.1(2)). The theorem now follows from Corollary 5.3.
\qed\enddemo 

Another interesting metric characterization involves extensions of
 geodesics. A local geodesic $\g$ ending at $x$ {\it extends
 discretely} if the set of directions in which $\g$ can be
 geodesically continued through $x$ is a non-empty discrete subset of
 $\lk (x,X)$, or equivalently, if the set of points $v \in \lk(x,X)$
 at distance $\geq \pi$ from $\g_x$ is non-empty and discrete. We say
 a geodesic metric space $X$ has the {\it discrete extension property}
 if $\g$ extends discretely for every local geodesic $\g$.

\proclaim{Theorem 7.5} Suppose $X$ is a connected, locally spherical
 (respectively Euclidean) space of dimension $n \geq 2$, and suppose
\roster
\item $X$ is CAT(1) (respectively CAT(0)),
\item $X$ has the discrete extension property.
\endroster
Then $X$ is isometric to a spherical building (respectively metric Euclidean
building).
\endproclaim

\demo{Proof} First note that if $X$ has the discrete extension
property, then so does $L_x = \lk(x,X)$ for every $x$.
For if $v \in L_x$ and $\g$ is a geodesic emanating from $x$ in
direction $v$, then for any point $y=\g(t)$ in a conelike neighborhood 
of $x$, $L_y = \s^0 *\lk(v,L_x)$. Thus, if every point in $L_y$ has a
non-empty discrete set of points at distance $\geq \pi$, then the same
holds in $\lk(v,L_x)$. (See the appendix of \cite{CD} for details on
distances in spherical suspensions.)

The discrete extension property also implies that $L_x$ is connected.
For if $\sigma$ is a spherical $(n-1)$--cell in $L_x$ and $v \in L_x$ is not
in the connected component of $\sigma$, then all of $\sigma$ has
distance $\geq \pi$ from $v$.

The theorem now follows by induction on $n$. If $n=2$, then the
discrete extension property implies that $L_x$ is connected, has
diameter $\pi$, and every vertex has valence at least 2. Thus, by
Theorem 6.1, and the remark following it, $L_x$ is isometric to a
1-dimensional spherical building for all $x$. By Theorem 3.8,
(respectively Corollary 5.3) $X$ is isometric to a spherical (respectively metric
Euclidean) building.

If $n > 2$, then $L_x$ is a connected, locally spherical space
satisfying conditions (1) and (2) of the theorem. By induction $L_x$
is isometric to a spherical building for every $x$, and the conclusion
follows from Theorem 3.8 (respectively Corollary 5.3).
\qed\enddemo


\Refs


\ref\key{Ba}
\by W Ballmann
\book Lectures on Spaces of Nonpositive Curvature
\bookinfo DMV Seminar 25,
\publ Birkh\"auser \yr 1995\endref

\ref\key{BB1}
\by W Ballmann \by M Brin
\paper Orbihedra of nonpositive curvature
\jour Inst. Hautes \'Et\-udes Sci. Publ. Math.
\vol 82 \yr 1996 \pages 169--209
\endref

\ref\key{BB2}
\by W Ballmann \by M Brin
\paper Diameter rigidity of spherical polyhedra
\jour Duke Math. J.
\vol 97 \yr 1999 \pages 235--259
\endref

\ref\key{BH}
\by M Bridson \by A Haefliger
\book Metric Spaces of Non-positive Curvature
\publ Springer--Verlag \yr 1999\endref

\ref\key{Bo}
\by N Bourbaki
\book Groupes et Algebr\`es de Lie
\bookinfo Ch. 4--6,
\publ Masson, \publaddr Paris \yr 1981\endref

\ref\key{Br}
\by K Brown
\book Buildings
\publ Springer--Verlag, \publaddr New York \yr 1989\endref

\ref\key{CD}
\by R Charney \by M\,W Davis
\paper Singular metrics of nonpositive curvature on bran\-ched covers 
of Riemannian manifolds
\jour Amer. J. Math. \vol 115 \yr 1993 \pages 929--1009
\endref

\ref\key{D}
\by M\,W Davis
\paper Buildings are CAT(0)
\inbook Geometry and Cohomology in Group Theory
\bookinfo ed. by P Kropholler, G Niblo and R St\"ohr, 
LMS Lecture Notes Series 252, 
\publ Cambridge University Press \yr1998\endref

\ref\key{G}
\by M Gromov
\paper Hyperbolic groups
\inbook Essays in Group Theory
\bookinfo ed. by S\,M Gersten, MSRI Publ. 8,
\publ  Springer--Verlag, \publaddr New York \yr 1987 \pages 75--264
\endref

\ref\key{K}
\by W\,M Kantor
\paper Generalized polygons, SCAB's and GAB's 
\inbook Buildings and the Geometry of Diagrams
\bookinfo LNM 1181, 
\publ Springer--Verlag, \publaddr New York \yr 1986 \pages 79--158
\endref

\ref\key{KL}
\by B Kleiner \by B Leeb
\paper Rigidity of quasi-isometries for symmetric spaces and Euclidean
buildings 
\jour Publ. IHES \vol 86 \yr 1997 \pages 115--197
\endref

\ref\key{N}
\by A Neumaier
\paper Some aporadic geometries related to PG(3,2)
\jour Arch. Math. \vol 42 \yr 1984 \pages 89--96
\endref

\ref\key{R}
\by M Ronan
\book Lectures on Buildings
\bookinfo Perspectives in Mathematics, Vol. 7,
\publ Academic Press, \publaddr San Diego \yr 1989\endref

\ref\key{T}
\by J Tits
\paper A local approach to buildings
\inbook The Geometric Vein, Coxeter Festschrift
\bookinfo ed. by C Davis, B Gr\"uenbaum and F\,A Sherk,
\publ  Springer--Verlag, \publaddr New York \yr 1981 \pages 519--547
\endref
\endRefs
\enddocument